\documentclass{article}
\usepackage[utf8]{inputenc}

\usepackage{amsfonts}
\usepackage{amsmath}
\usepackage{amssymb}
\usepackage{mathtools}
\usepackage{xcolor}
\usepackage{fullpage} 
\usepackage{hyperref}
\usepackage{cleveref}
\usepackage{graphicx}
\usepackage{caption}
\usepackage{subcaption}
\usepackage{authblk}
\usepackage{soul}
\usepackage{orcidlink}

\def\sgn{\operatorname{sgn}}

\newtheorem{theo}{Theorem}

\newtheorem{deff}{Definition}

\newtheorem{lemma}{Lemma}

\newtheorem{corollary}{Corollary}

\newtheorem{theorem}{Theorem}

\newenvironment{proof}{\paragraph{Proof:}}{\hfill$\square$}

\newcommand{\setword}[2]{%
  \phantomsection
  #1\def\@currentlabel{\unexpanded{#1}}\label{#2}%
}

\makeatletter

\def\blfootnote{\gdef\@thefnmark{}\@footnotetext}

\makeatother

\newcommand{\argmin}{\operatornamewithlimits{argmin}}
\newcommand{\diag}{\operatorname{diag}}

\usepackage[backend=biber,style=vancouver]{biblatex}
\begin{document}

    \title{Auto-Calibration and Biconvex Compressive Sensing with Applications to Parallel MRI
    }
    
    \author{%
    Yuan Ni \,\orcidlink{0000-0002-5797-5234} \thanks{Department of Mathematics, University of California Davis, Davis, CA 95616-5270, USA \texttt{yuani@ucdavis.edu}} %
    \, and   
    Thomas Strohmer \,\orcidlink{0000-0003-2029-3317}\thanks{Department of Mathematics
and Center of Data Science and Artificial Intelligence Research, University of California Davis,
Davis, CA 95616-5270, USA \texttt{strohmer@math.ucdavis.edu}}%
    }
  
    \date{}
    
    \maketitle
 
\newcommand{\bb}[1]{\textbf{#1}}
\newcommand{\dd}{\text{diag}}

    
   \begin{abstract}
       We study an auto-calibration problem in which a transform-sparse signal is acquired via compressive sensing by multiple sensors in parallel, but with unknown calibration parameters of the sensors. This inverse problem has an important application in pMRI reconstruction, where the calibration parameters of the receiver coils are often difficult and costly to obtain explicitly, but nonetheless are a fundamental requirement for high-precision reconstructions. Most auto-calibration strategies for this problem involve solving a challenging biconvex optimization problem, which lacks reconstruction guarantees. In this work, we transform the auto-calibrated parallel compressive sensing problem to a convex optimization problem using the idea of `lifting'. By exploiting sparsity structures in the signal and the redundancy introduced by multiple sensors, we solve a mixed-norm minimization problem to recover the underlying signal and the sensing parameters simultaneously. Our method provides robust and stable recovery guarantees that take into account the presence of noise and sparsity deficiencies in the signals. As such, it offers a theoretically guaranteed approach to auto-calibrated parallel imaging in MRI under appropriate assumptions. Applications in compressive sensing pMRI are discussed, and numerical experiments using real and simulated MRI data are presented to support our theoretical results. 
       \\
       \\
       {\bf Keywords:}
        Self-calibration, Compressive sensing, Convex optimization, Random matrices, Parallel MRI (pMRI), Inverse problems.
   \end{abstract}

                
   \section{Introduction}
        We frequently encounter challenges with imperfect sensing, where accurately calibrated sensors are critical for high-precision measurements and reconstructions. For many applications, explicit calibration is often difficult and expensive to carry out in practice, posing a major roadblock in many scientific and technological endeavors where accurate and reliable measurements are essential. The idea of auto-calibration (calibration-less/self-calibration) avoids the difficulties and inaccuracy associated with explicit estimations by equipping the sensors and systems with the capability to automatically derive the sensing information from its collected signals while performing the intended system function. Auto-calibration is a field of research that finds successful applications across diverse fields such as wireless communication, computer vision, remote sensing, biomedical imaging, and more. This field encompasses a wide range of techniques, including blind deconvolution, blind phase retrieval, blind rotation estimation, and direction-of-arrival estimation, among others. 
\\

\noindent
The inverse problem for auto-calibration is generally difficult to solve. In practice, most self-calibration algorithms rely on joint estimation of both the signals of interest and the calibration parameters using standard techniques such as maximum likelihood estimation, alternating minimization, or convex/nonconvex optimization. In this work, we extend the line of work that combines convex optimization and compression sensing in solving various auto-calibration tasks \cite{ling_strohmer_2015,Candès_Eldar_Strohmer_Voroninski_2013,Ahmed_Recht_Romberg_2014,flinth_2017}. A common scheme in such approaches includes first transforming the biconvex constraints to convex constraints using the concept of `lifting' \cite{Candès_Eldar_Strohmer_Voroninski_2013}; then, drawing on ideas from compressive sensing to exploit redundancy/sparsity priors in the signal model and/or the calibration model while being consistent with the convex constraints.
\\

\noindent
In this work, we consider a type of auto-calibration problem where the same underlying signal is sensed by multiple sensors in parallel with unknown calibration parameters. Specifically, we are concerned with the following problem:
\begin{equation}\label{eq:1}
    \mathbf{y}_i = \mathbf{F}_\Omega (\mathbf{s}_i \odot \mathbf{x}) + \boldsymbol{\omega}_i, i = 1 \cdots C ,
\end{equation}
where $\mathbf{y}_i \in \mathbb{C}^{L \times 1}$ are the measurements, $\mathbf{F}_\Omega \in \mathbb{C}^{L \times N} \text{ with } L\leq N$ is the partial Fourier matrix with rows belonging to a subset $\Omega \subset [N]$, $\mathbf{s}_i \in \mathbb{C}^{N\times 1}$ is the vector of unknown calibration parameters, $\mathbf{x}\in \mathbb{C}^{N\times 1}$ is our signal of interest, $\boldsymbol{\omega}_i \in \mathbb{C}^{L\times 1}$ is the additive noise. We have $C$ sets of measurements. Here, $\odot$ is the point-wise multiplication. Our goal is to simultaneously reconstruct the signal $\bb{x}$ and the parameters related to the sensor profiles (calibration parameters) $\bb{h}_i$. 
\\
\\
\noindent
Such a forward model appears, for instance, in blind parallel MRI (pMRI). Specifically, the $i$-th coil measurement $\mathbf{y}_i$ is a noisy sample of the spatial Fourier transform of an object $\mathbf{x}$, weighted by an unknown coil sensitivity $\mathbf{s}_i$ that corresponds to the $i$-th sensing coil. In parallel imaging, $C$ coils are used to take measurements simultaneously for the same underlying object $\mathbf{x}$. The joint estimation of coil sensitivities and the object of interest leads to the nonlinear inverse problem of the form \eqref{eq:1}.
\\

\noindent
Before proceeding with further analyzing the model, we briefly discuss the identifiability and dimension of the problem. First, if $\mathbf{s}_i$ and $\mathbf{x}$ form a pair of solutions to \eqref{eq:1}, then for any complex vector $\mathbf{v}$ with nonzero entries, $\mathbf{s}_i \odot \mathbf{v}$ and $\mathbf{x} \odot 1/\mathbf{v}$ are a pair of solutions. Furthermore, with $L \leq N$, \eqref{eq:1} has $CL$ number of measurements and $(C+1)N$ number of unknowns. Thus, it is impossible to recover all $\mathbf{s}_i$ and $\mathbf{x}$ without making further assumptions. 
\\

\noindent
Of particular interest is the case when $\mathbf{x}$ is {\em transform-sparse}, a characteristic often observed in real-world signals \footnote{This notion can be easily generalized to {\em approximately} transform-sparse.}. That is, there is some known sparsifying transformation $\boldsymbol{\Psi} \in \mathbb{C}^{N\times N}$ such that $\boldsymbol{\Psi}\bb{x}$ is n-sparse with $n < N$. Common choices for $\boldsymbol{\Psi}$ include the well-known DCT basis and the wavelet basis. Despite the sparsity assumption, the available measurements are still fewer than the number of unknowns (i.e., $CL \leq CN+n$). We impose an additional subspace assumption, $\mathbf{s}_i = \mathbf{B}\mathbf{h}_i$, which can be interpreted as that the calibration parameters all lie in some known subspace spanned by a tall matrix $\mathbf{B}$. Practical choices for such $\mathbf{B}$ can be a low-order polynomial basis or a low-frequency sinusoidal basis. Such choices promote smoothness in the sensing parameters \cite{pruessmann_weiger_scheidegger_boesiger_1999,Ying_Sheng_2007,Guerquin-Kern_Lejeune_Pruessmann_Unser_2012,Uecker_Hohage_Block_Frahm_2008}. Altogether, the forward model becomes:
\begin{equation}\label{eq:2}
    \mathbf{y_i} = \mathbf{F}_\Omega (\mathbf{s_i} \odot \mathbf{x}) + \boldsymbol{\omega}_i, \mathbf{s_i} = \mathbf{B}\mathbf{h_i}, \mathbf{\Psi}\mathbf{x} = \mathbf{z}, \mathbf{z} \text{ is sparse }, i = 1 \cdots C .
\end{equation}
Here, $\mathbf{B}\in \mathbb{C}^{N \times k}, k < N$, $\mathbf{\Psi} \in \mathbb{C}^{N \times N}$ are known transformations. 
\\

\noindent
Although the above model may appear straightforward and specific, providing a rigorous analysis is nontrivial. The primary challenge arises from carefully analyzing the incoherence of the forward matrix, which involves the interactions between $\mathbf{F}_\Omega$, $\mathbf{B}$, and $\boldsymbol{\Psi}$. Our goal is to understand the requirements on the number of samples $L$ and the number of channels $C$ for robust and stable recovery of all $\mathbf{s}_i$ and $\mathbf{x}$. We will discuss particular choices of $\mathbf{F}_\Omega$, $\mathbf{B}$, and $\boldsymbol{\Psi}$ that align with the pMRI application in detail in Section ~\ref{2}.
\\

\noindent
An important observation is that Parallel Imaging (PI) studies can sometimes appear ``noisy". This is largely because the primary objective of PI techniques is to reduce imaging time, which leads to acquiring fewer data points per coil due to the need to reduce imaging time. Consequently, this reduction results in a lower signal-to-noise ratio (SNR). As a result, the SNR of PI sequences is typically lower than that of comparable non-PI sequences \cite{Aja-Fernández_Vegas-Sánchez-Ferrero_Tristán-Vega_2014,Robson_Grant_Madhuranthakam_Lattanzi_Sodickson_McKenzie_2008}.
\\

\noindent
The SNR of PI sequences is affected by the same factors that influence non-PI sequences, such as field strength, magnet hardware, coil loading, the type of tissue being imaged, pulse sequence, timing parameters, voxel volume, the number of phase-encoding steps, receiver bandwidth, and the number of signals averaged. Additionally, the SNR for PI sequences ($\text{SNR}_{\text{parallel}}$) is further reduced by two extra factors, as shown in the equation below:
$$
\text{SNR}_{\text{parallel}} = \text{SNR}_{\text{non-PI}} / (\sqrt{R} \cdot g),$$
where $R$ denotes the reduction or acceleration factor for undersampling, and $g$ represents a spatially dependent term known as the geometry factor. The term $g$ relates to the number, size, and orientation of the surface coil elements and can be considered a measure of coil separation \cite{Aja-Fernández_Vegas-Sánchez-Ferrero_Tristán-Vega_2014,Robson_Grant_Madhuranthakam_Lattanzi_Sodickson_McKenzie_2008}.

\subsection{Notation and outline}

We denote vectors and matrices by bold font letters (e.g. \bb{x} or \bb{X}) and scalars by regular font or Greek symbols. We denote the circular convolution as $*$ and the Kronecker product as $\otimes$. For any integer $N$, we denote the set $\{1,\cdots,N\}$ as $[N]$. For any vector $\bb{x} \in \mathbb{C}^N$ and set $S \subset [N]$, define $\bb{x}_S = \mathcal{P}_S\bb{x} \in \mathbb{C}^N$ the orthogonal projection of $\bb{x}$ onto set $S$. (e.g. $\bb{x}_S(i)= \bb{x}(i),i \in S$ and $\bb{x}_S(i)= 0, \text{for }i \in S^c$). Similarly, define the orthogonal projection of matrix $\bb{X}\in \mathbb{C}^{N_1 \times N_2}$ to set $T \subset [N_1] \times [N_2] $ by $\bb{X}_T = \mathcal{P}_T\bb{X} \in \mathbb{C}^{N_1 \times N_2}$, such that $\bb{X}_T(i,j)= \bb{X}(i,j),(i,j) \in T$ and $\bb{X}_T(i,j)= 0,\text{otherwise}$. We denote the N by N identity matrix as $\bb{I}_N$. For any complex vector $\bb{x}$ or matrix $\bb{X}$, $\bar{\bb{x}}$ or $\bar{\bb{X}}$ means their conjugate, respectively. $vec(\bb{X})$ means the vectorization of a matrix $\bb{X}$ in the column-wise order into a vector. For any set $S$, we denote $S^c$ as its complementary set and $|S| = \text{Card}(S)$. Informally, we refer to a matrix $X \in \mathbb{C}^{N_1\times N_2}$ being incoherent if the 2-norm of its rows does not vary too much.
\\

\noindent
The paper is organized as follows: In Section \ref{2}, we briefly discuss our contributions in relation to the state of the art in auto-calibration and compressive sensing pMRI. In Section \ref{3}, we introduce the problem setup. We present our main results in Section~\ref{4} and provide the details of the proofs in Section~\ref{6}. Numerical experiments are shown in Section~\ref{5}.
   \section{Related work}\label{2}            
        The work is inspired by papers related to bringing auto-calibration problems into the framework of biconvex compressive sensing. Many auto-calibration problems can be seen as jointly recovering the signal and the calibration parameters from bilinear measurements. Specifically, they can be seen as variations of recovering two unknown signals $\mathbf{x}$ and $\mathbf{z}$ from bilinear measurements,
\begin{equation}\label{eq:3}
    y_i = \langle \mathbf{x}, \mathbf{b_i}\rangle \langle \mathbf{z},\mathbf{\bar{a}_i} \rangle, i = 1 \cdots m,
\end{equation}
where $\mathbf{x}, \mathbf{b}_i \in \mathbb{C}^{N_1\times 1} , \mathbf{z},\mathbf{a}_i \in \mathbb{C}^{N_2\times 1}$. 
\\

\noindent
Various auto-calibration problems of similar setups have been explored, including blind deconvolution in \cite{Ahmed_Recht_Romberg_2014}, direction-of-arrival estimation as seen in \cite{ling_strohmer_2015}, blind deconvolution and demixing \cite{flinth_2017}, phase retrieval via matrix completion \cite{Candès_Eldar_Strohmer_Voroninski_2013}, among others. The analysis of the auto-calibration problems generally differs in the following two aspects:
\begin{itemize}
    \item Additional assumptions must be made about the signal, the calibration parameters, or both to address the ill-posed nature of the problem. Popular assumptions may vary in several ways, either regarding the signal, the sensing parameters, or both: whether they assume sparsity individually or jointly, sparsity in the natural parameter space or after transformations, subspace constraints, or restrictions imposed by sign constraints.
    \item Specific application cases impose different forward models due to physical or other constraints. For example, the choice of sensing matrices - whether sub-Gaussian (e.g., $b_i$ or $a_i$ are Gaussian), structured random (e.g., taking random samples from the Fourier transform) or deterministic.
\end{itemize}

\noindent
Two primary approaches to solving the bilinear optimization problem are as follows: one involves convexifying it through methods such as linearization or lifting \cite{Ahmed_Recht_Romberg_2014,ling_strohmer_2015,Ling_Strohmer_2018,flinth_2017,Candès_Eldar_Strohmer_Voroninski_2013}. The other involves optimization in the natural parameter space using alternating minimization and gradient descent algorithms \cite{Lee_Tian_Romberg_2018,Li_Ling_Strohmer_Wei_2019,Qu_Li_Zhu_2019}. 
\\

\noindent
For a general discussion on the issue of injectivity and the principle of identifiability for bi-linear compressive sensing, one may refer to \cite{Kech_Krahmer_2017,Choudhary_Mitra_2013}. 
\subsection{Related work in auto-calibration and compressive sensing}
It seems natural and tempting to recover $(\bb{h}_1,\cdots,\bb{h}_C, \bb{z})$ obeying \eqref{eq:2} by  solving the optimization problem
\begin{equation}\label{nonconvex}
\min\limits_{\bb{h}_1,\cdots,\bb{h}_C, \bb{z}} \sum\limits_{i} \|\bb{y}_i - \bb{F}_{\Omega}(\bb{B}\bb{h}_i \odot \boldsymbol{\Psi}^*\bb{z})\|_2^2.
\end{equation}
\noindent
However, the above nonlinear least squares problem is a challenging optimization problem since it is highly nonconvex and most of the available algorithms, such as alternating minimization and gradient descent, may suffer from easily getting trapped in some local minima. 
\\

\noindent
A popular convex optimization approach involves transforming the bi-convex problem \eqref{eq:3} to the convex problem using the idea of `lifting'. Define the matrix $\mathbf{X} := \mathbf{x}\mathbf{z}^T$, we can write the bilinear equations as a linear equation with respect to $\mathbf{X}$ using some linear algebra,
\begin{equation}
    \langle \mathbf{x}, \mathbf{b_i} \rangle \langle \mathbf{z},\mathbf{\bar{a}_i} \rangle = \mathbf{b}_i^* \mathbf{x} \mathbf{z}^T \mathbf{a}_i = \mathbf{b}_i^* \mathbf{X} \mathbf{a}_i = (\mathbf{a}_i^T \otimes \mathbf{b}_i^*) \operatorname {vec}(\mathbf{X}).
\end{equation}
\noindent
Hence, by defining the linear operator $\mathcal{A} :\mathbb{C}^{N_1 \times N_2}\rightarrow \mathbb{C}^{m}$, the system of equations becomes
\begin{equation}
[y_1,\cdots, y_m] = \mathcal{A}(\mathbf{X}) := \{ \mathbf{b}_i^* \mathbf{X} \mathbf{a}_i \}_{i=1}^m .
\end{equation}
The original bilinear equations of dimension $(N_1 + N_2)$ are `lifted' to linear equations of underlying dimension  $N_1 \times N_2$.
\\

\noindent
An important auto-calibration problem pertaining to sparse signals arises from the task of image reconstruction using randomly coded masks \cite{Bahmani_Romberg_2015,Tang_Recht_2014}. The forward model for this auto-calibration problem is given by: 
\begin{equation}\label{ling}
    \mathbf{y} = \text{diag}(\mathbf{B}\mathbf{h})\mathbf{A}\mathbf{x} + \boldsymbol{\omega},  
\end{equation}
where $\mathbf{h}\in \mathbb{C}^{k}$, $\mathbf{y}\in \mathbb{C}^{L}$, $\mathbf{A}\in \mathbb{C}^{L \times N}, L \leq N$ and $\|\mathbf{x}\|_0 = n$. In this model, the matrix $\bb{A}$ acts as the sensing matrix. The differences between the forward model \eqref{ling} and the current work \eqref{eq:1} are as follows: 1. The blind sensing (i.e., $\text{diag}(\mathbf{Bh})$) is applied after compressed sensed signals (i.e., $\mathbf{Ax}$). 2. The signal is sparse in its natural domain rather than being transform-sparse. 
\\

\noindent
The SparseLift framework~\cite{ling_strohmer_2015} solves the auto-calibration problem in \eqref{ling} by combining `lifting' and bi-convex compressive sensing. The authors in \cite{ling_strohmer_2015} first transform the biconvex constraint \eqref{ling} to a linear constraint with respect to $\bb{X} = \bb{h}\bb{x}^T$. Exploiting the sparsity in $\mathbf{x}$ and consequently within $\bb{X}$, the recovery solves a convex relaxation by minimizing the $\|\mathbf{X}\|_1$ subject to the linear constraints. The authors show the recovery of $\mathbf{X}_0 = \mathbf{h}_0\mathbf{x}_0^T$ is guaranteed with high probability when the sensing matrix $\mathbf{A}$ is either a Gaussian random matrix or with its rows chosen uniformly at random with replacement from the discrete Fourier transformation (DFT) matrix, when the number of measurements $L$ is of the order of $\mathcal{O}(kn\log^2(L))$. Note that the case where the signal is only approximately sparse is not analyzed. In a related work, the author in \cite{flinth_2017} extends the results of SparseLift to the case where one observes a summation of such measurements from multiple sparse signals, $\mathbf{y} = \sum\limits_{i\in [r] }\text{diag}(\mathbf{B}\mathbf{h}_i)\mathbf{A}\mathbf{x}_i + \boldsymbol{\omega}$. This extension improves recovery guarantees by solving a $\ell_{1,2}$ minimization by promoting block sparsity in the lifted signals. 
\\

\noindent
The idea of `lifting' also applies to another class of self-calibration problems, notably the blind deconvolution problem \cite{Ahmed_Recht_Romberg_2014,Levin_Weiss_Durand_Freeman_2011}. In the general form of blind deconvolution, one measures the convolution of a signal with an unknown filter. That is, $\bb{y} = \mathbf{I}_{\Omega}(\bb{s} * \bb{x}) + \bb{n}$ where both $\bb{s}$ and $\bb{x}$ are unknown. Notice that the forward model for our problem in \eqref{eq:2} can be written in the form of a convolution $ \mathbf{y} = \mathbf{I}_{\Omega}(\mathbf{F}\mathbf{Bh} * \mathbf{F}\boldsymbol{\Psi^*}\mathbf{z}) + \bb{n}$ using the convolution theorem. When the Fourier space is fully sampled (i.e. $\mathbf{I}_{\Omega} = \mathbf{I}$ ), the problem is equivalent to recovery from 
\begin{equation} \label{eq:4}
\mathbf{F}^{-1}\mathbf{y} = \mathbf{Bh} * \boldsymbol{\Psi^*}\mathbf{z} + \mathbf{w},
\end{equation} where $\mathbf{w}$ is also Gaussian. The work by Ahmed et al. \cite{Ahmed_Recht_Romberg_2014} provides a theoretically guaranteed solution to the fully sampled case in \eqref{eq:4}. Using lifting, the authors solve a nuclear norm minimization problem by exploiting low-rankness in the lifted signal. If the active coefficients of $\bb{z}$ are known and $\|\bb{z}\|_0 = n$, the main theorem of \cite{Ahmed_Recht_Romberg_2014} guarantees a stable recovery when the number of measurements is, up to log factor, of the order $\mathcal{O}(n+k)$ when $\boldsymbol{\Psi}$ is a random Gaussian matrix and $\bb{B}$ is incoherent. For the more difficult case when the convolution is only subsampled, the authors in \cite{Lee_Li_Junge_Bresler_2017} provide rigorous analysis when both $\bb{s}$ and $\bb{x}$ are transform sparse (i.e., there exists $\boldsymbol{\Phi}$ and $\boldsymbol{\Theta}$ such that $\boldsymbol{\Phi}\bb{x}$ and $\boldsymbol{\Theta}\bb{s}$ are sparse). It provides an iterative algorithm for robust recovery guarantees with near optimal sample complexity. However, the recovery can only be guaranteed when both the sparsifying transformations (i.e. $\boldsymbol{\Phi}$ and $\boldsymbol{\Theta}$) are random Gaussians. As seen in many theories, subgaussian random matrices provide optimal measurement matrices for compressive sensing \cite{foucart_rauhut_2015}. However, note that for our application case, $\boldsymbol{\Phi} = \boldsymbol{\Psi}\bb{F}^{-1}$ and $\boldsymbol{\Theta} = \bb{B}^{*}\bb{F}^{-1}$, both impose structures due to physical constraints.
\\

\noindent
Similar to the parallel imaging (i.e., $C > 1$ in \eqref{eq:2}), from a practical viewpoint, a relevant scenario focuses on self-calibration from multiple snapshots. For example, blind image restoration from multiple filters \cite{Harikumar_Bresler_1999} and self-calibration model for sensors \cite{Balzano_Nowak_2007,Ling_Strohmer_2018,Gribonval_Chardon_Daudet_2012,Bilen_Puy_Gribonval_Daudet_2014}. The goal here is to recover the unknown gains/phases $\mathbf{S} = \text{diag}(\mathbf{s})$ and a signal matrix $\mathbf{X} = [\mathbf{x}_1,\cdots,\mathbf{x}_p] $ from the measurement matrix $\mathbf{Y}=[\mathbf{y}_1,\cdots,\mathbf{y}_p]$ and \begin{equation}\label{eq:5}
     \mathbf{Y}= \mathbf{S}\mathbf{A}\mathbf{X}.
\end{equation} For this model, the sensing matrix $\mathbf{A}$ is fixed through the sensing process, just as in the pMRI case. In general, the method does not require sparse priors on the signal or the calibrations. Among this line of research, an interesting and somewhat related method involves using the idea of `linearization' to convexify the biconvex measurements \eqref{eq:5} into a linear one. One crucial assumption for the `linearization' scheme is that $\mathbf{s}$ does not have zero entries. Define $\mathbf{s} = \mathbf{1/d}$, the equations $ \mathbf{y}_l = \text{diag}(\mathbf{s})\mathbf{A}\mathbf{x_l}$ can be written as
\begin{equation}
    \mathbf{y}_l  \mathbf{d}= \mathbf{A}\mathbf{x_l},
\end{equation}
which is a linear equation in the concatenated vector $\begin{bmatrix}
\mathbf{d} \\ \mathbf{x}_l
\end{bmatrix}$. Compared to `lifting', an advantage of `linearization' is that the underlying problem dimension remains unchanged. The same idea of linearization can also be found in \cite{Gribonval_Chardon_Daudet_2012,Bilen_Puy_Gribonval_Daudet_2014,Balzano_Nowak_2007,Wang_Chi_2016}. In \cite{Ling_Strohmer_2018}, the authors use `linearization' and solve a least-square minimization problem for \eqref{eq:5}. Furthermore, theoretical recovery guarantee is provided when $\mathbf{A}$ is Gaussian. 
\\

\noindent
As observant readers already may have noticed, the `true dimension' of the problem of recovering $(\bb{h}_1,\cdots,\bb{h}_c,\bb{z})$  is $(kC + n)$. To the best of the authors' knowledge, currently there are no convex optimization schemes which succeed with such few measurements. In this context, the work by Oymak et al. \cite{Oymak_Jalali_Fazel_Eldar_Hassibi_2015} should be mentioned. The authors demonstrated nonconvex procedures (e.g., minimizing a linear combination of $\ell_0$, $\ell_{0,2}$ or rank of the lifted matrix) that succeed with high probability, when the number of measurements is (up to log-factors) slightly larger than the information-theoretical limits. However, the theories apply only to measurement ensembles with specific properties, such as matrices with sub-Gaussian rows, subsampled standard basis (e.g. in matrix completion), or quadratic measurements for phase retrieval. Also, these procedures are not efficiently computable. For other structured measurement ensembles, it is worth noting the works by Qu et al. and Ling et al. \cite{Qu_Li_Zhu_2019,Li_Ling_Strohmer_Wei_2019}. These studies address certain blind deconvolution problems with efficient algorithms and theoretical guarantees when the number of measurements is (up to log-factors) slightly larger than the information-theoretical limits. However, to the best of our knowledge, none of these methods directly apply to our measurement ensemble or require additional assumptions that do not align with our application case.

\subsection{Developments in compressive sensing methods in pMRI}\label{2.2}
In the literature on pMRI, combining compressive sensing and parallel imaging has gained significant importance due to the need for fast acquisition. Accelerated imaging with fewer samples using multiple receiver coils has been an active field of research since its beginning. In this section, we briefly discuss the development of compressive sensing approaches in pMRI that exploit transform sparsity in the time domain, with a focus on those most closely related to our approach.

\subsubsection{Relevant problem formulations and the optimization problems}\label{previous_methods}
In parallel imaging, the data are acquired from multiple channels simultaneously, but each at a rate lower than the Nyquist rate such that the data acquisition time is reduced. If only $1/R$ of the Fourier space is sampled per channel, the acquisition time is reduced by a factor of $R$. 
\\

\noindent
The measurement model for pMRI using $C$ coils/channels is the following,
\begin{equation}\label{pMRI}
    \bb{Ex} = \bb{y},
\end{equation}
where the encoding matrix $\bb{E}$ consists of the product of partial Fourier encoding with channel-specific sensitivity encoding over the image. Specifically,  $\bb{E} = \begin{bmatrix}
    \bb{F}_\Omega \bb{D}_1 
    \\
    \vdots
    \\
    \bb{F}_\Omega \bb{D}_C
\end{bmatrix}$ and $\bb{D}_i = \text{diag}(\bb{s}_i)$ where $\bb{s}_i$ is the coil sensitivity for the $i$-th coil. 
\\

\noindent
According to linear algebra and the generalized sampling theorem, the maximum reduction factor is equal to the number of channels if $\bb{E}$ is of full rank. Because these conditions are rarely guaranteed in practice, the possible reduction factor is much lower. In practice, the reduction factor is typically chosen between 2 and 6 \cite{Deshmane_Gulani_Griswold_Seiberlich_2012, Glockner_Hu_Stanley_Angelos_King_2005} for parallel imaging when the number of coils is large. Thus, at the boundary case, one would expect a threshold for the maximum reduction factor and the minimum number of channels necessary to ensure the invertibility of such a linear problem.
\\

\noindent
In addition to reducing sampling time through sensitivity encoding, further reduction in the number of samples can be achieved by exploiting redundancy in the signals. In this work, we specifically discuss time-domain sparsity, as magnetic resonance (MR) images typically exhibit sparse representation in some appropriate transform domain. Suppose it is assumed that the signal is \( n \)-sparse after some known transformation. One approach is to solve the following \( \ell_0 \) norm constrained least squares problem to directly exploit the \( n \)-sparsity of the signal:

\begin{equation}\label{l0}
    \min\limits_{\bb{x}} \|\bb{F}_{\Omega}\bb{x} - \bb{y}\|_2^2 \quad \text{subject to } 
    \|\boldsymbol{\Psi}\bb{x}\|_0 = n.
\end{equation}
As a convex relaxation to \( \ell_0 \) minimization, the following \( \ell_1 \) minimization is often utilized for practical compressive sensing MRI applications. Since Lustig et al.~\cite{lustig_donoho_pauly_2007} first demonstrated the compressive sensing approach for MRI reconstruction, it has become one of the most important tools used for modern MR imaging research. The original $\ell_1$ norm minimization problem solves the following convex program,
\begin{equation}\label{sparseMRI}
    \min\limits_{\bb{x}} \|\boldsymbol{\Psi}\bb{x}\|_1 \quad \text{subject to } \bb{F}_{\Omega}\bb{x} = \bb{y},
\end{equation}
where the signal $\bb{x}$ is assumed to be transform sparse. In the MRI literature, it has been demonstrated that the discrete cosine transformation, the wavelet transform, and the finite difference transform are the most commonly utilized choices for sparsifying MR images \cite{shin_larson_ohliger_elad_pauly_vigneron_lustig_2014,Ravishankar_Bresler_2015,Chun_Adcock_Talavage_2014}. The choice of such $\boldsymbol{\Psi}$ is related to the field of research in sparsifying transformations \cite{Elad_Milanfar_Rubinstein_2007,Pratt_Kane_Andrews_1969}. More recently, attention has turned to learning the sparsifying transformations from data \cite{Aharon_Elad_Bruckstein_2006,Ravishankar_Bresler_2012}.
\\

\noindent
In terms of combining CS with parallel imaging, the first applications were demonstrated in~\cite{Liang_Liu_Wang_Ying_2009,Bo_Liu_Yi_Ming_Zou_Leslie_Ying_2008,Ji_Chen_Zhao_Tao_Lang_2008,Wu_Millane_Watts_Bones_2010}. They can be seen as a direct extension of \eqref{sparseMRI} to the parallel case,
\begin{equation}\label{sparseSENSE}
     \min\limits_{\bb{x}} \|\boldsymbol{\Psi}\bb{x}\|_1 +\alpha \|\bb{x}\|_{\text{TV}} \quad \text{subject to } \bb{E} \bb{x} = \bb{y},
\end{equation}
where $\bb{D}_i = \text{diag}(\bb{s}_i)$ is the diagonal matrix of known calibration parameters of the $i$-th coil.
\\

\noindent
Note that both the convex relaxations in \eqref{sparseMRI} \eqref{sparseSENSE} and nonconvex problem in \eqref{l0} assume that the coil sensitivities are known. In practice, coil sensitivities are typically determined either through a separate pre-scan \cite{pruessmann_weiger_scheidegger_boesiger_1999,mathematics_and_physics_of_emerging_biomedical_imaging_1996,Liang_Lauterbur_2000} or by interpolating from a fully sampled region of k-space (Fourier space) \cite{griswold_jakob_heidemann_nittka_jellus_wang_kiefer_haase_2002,McKenzie_Yeh_Ohliger_Price_Sodickson_2002,Morrison_Jacob_Do_2007}. However, both methods have inherent drawbacks. The pre-scan approach is susceptible to motion artifacts, environmental fluctuations, and requires additional acquisition time. The interpolation method, while faster, may introduce errors due to imperfect k-space sampling or interpolation techniques. Moreover, both approaches may necessitate the use of filtering and kernel selection to obtain reasonable estimates of the coil sensitivities, which can lead to variations in the resulting reconstructions.
\\

\noindent
Given these challenges, auto-calibrated reconstruction methods have gained attention. In particular, it considers recovering all $\bb{s}_i$ and $\bb{x}$ from the measurement model in \eqref{pMRI}. However, the linear system is highly underdetermined and ill-posed in this case. Indeed, if $\mathbf{s}_i$ and $\mathbf{x}$ form a pair of solutions to \eqref{pMRI}, then for any complex vector $\mathbf{v}$ with nonzero entries, $\mathbf{s}_i \odot \mathbf{v}$ and $\mathbf{x} \odot 1/\mathbf{v}$ form a pair of solutions. 
\\

\noindent
Prior knowledge about the object and the coil sensitivities must be incorporated to render reasonable solutions. As formulated in Joint Image Reconstruction and Sensitivity Estimation in SENSE (JSENSE) \cite{Ying_Sheng_2007} and Regularized Nonlinear Inversion (NLINV) \cite{Uecker_Hohage_Block_Frahm_2008}, the smoothness in the coil sensitivities is exploited. In JSENSE, the coil sensitivities are assumed to be expanded by a low-order polynomial basis (i.e. $\bb{B}$) and the following bi-convex optimization problem is solved via alternating minimization,
\begin{equation}\label{JSENSE}
    \arg\min \limits_{\mathbf{x},\mathbf{h_j}} \Sigma_{j=1}^{C}\|\mathbf{y}-\mathcal{F}_{\Omega}\{ \bb{B}\mathbf{h}_j \odot \mathbf{x} \}\|_2^2.
\end{equation}
In the NLINV formulation, they solve the following bi-convex optimization problem iteratively using a Gauss Newton type method,
\begin{equation}\label{NLIVE}
    \arg\min \limits_{\mathbf{x},\mathbf{s_j}} \Sigma_{j=1}^{C}\|\mathbf{y}-\mathcal{F}_{\Omega}\{ \mathbf{s}_j \odot \mathbf{x} \}\|_2^2 + \alpha (\Sigma_{j=1}^{C} \|\mathbf{W}\mathbf{s}_j\|_2^2+ \|\mathbf{x}\|_2^2),
\end{equation}
where $\mathbf{W} \in \mathbb{C}^{N\times N}$ is an invertible preconditioning matrix that penalizes high frequencies, thereby enforcing smoothness in the coil sensitivities. Other methods to enforce smoothness in the coil include expanding the coil sensitivities using low-order polynomial or sinusoidal basis functions, a method first introduced in the original SENSE paper \cite{pruessmann_weiger_scheidegger_boesiger_1999} and later adopted in the analytical phantom package \cite{Guerquin-Kern_Lejeune_Pruessmann_Unser_2012}.
\\

\noindent
An extension of the NLINV \eqref{NLIVE} framework is the ENLIVE \cite{Holme_Rosenzweig_Ong_Wilke_Lustig_Uecker_2019}, which solves the following relaxed joint nonconvex optimization problem over $k$ sets of measurements,
\begin{equation}\label{ENLIVE}
    \arg\min \limits_{\mathbf{x^i},\mathbf{s_j^i}} \Sigma_{j=1}^{C}\|\mathbf{y}_i-\mathcal{F}_{\Omega}\{ \Sigma_{i=1}^k \mathbf{s}_j^i \odot \mathbf{x}^i \}\|_2^2 + \alpha \Sigma_{i=1}^k (\Sigma_{j=1}^{C} \|\mathbf{W}\mathbf{s}_j^i\|_2^2+ \|\mathbf{x}^i\|_2^2).
\end{equation}

\noindent
Here, k sets of images $\mathbf{x}^i$ are used instead of one single image $\mathbf{x}$ as in \eqref{NLIVE} to account for model violations. In a post-processing step, an average image using the $k$ sets of solutions is calculated as the solution.
\\

\noindent
Of particular interest to us is that ENLIVE can be related to a nuclear norm minimization problem that exploits the low-rankness in a `lifted' domain. Specifically, by defining the matrices 
$$\mathbf{y} = \begin{bmatrix}
    \mathbf{y}_1 & \cdots & \mathbf{y}_C
\end{bmatrix}, \mathbf{U} = \begin{bmatrix} \mathbf{x}^1 & \cdots \mathbf{x}^k
\end{bmatrix} \in \mathbb{C}^{N\times k}, \mathbf{V} = \begin{bmatrix} 
\mathbf{W}\mathbf{s}^1_1 & \cdots \mathbf{W} \mathbf{s}^k_1 \\
\vdots & \vdots & 
\\
\mathbf{W}\mathbf{s}^1_C & \cdots \mathbf{W} \mathbf{s}^k_C
\end{bmatrix} \in \mathbb{C}^{CN \times k},$$ the bi-linear measurements can be transformed into linear measurements in terms of a rank-k matrix $\bb{U}\bb{V}^T$:

\begin{equation}
    \arg\min \limits_{\mathbf{U},\mathbf{V}} \|\mathbf{y}-\mathcal{A}(\mathbf{U}\mathbf{V}^T)\|_2^2 + \alpha (\|\mathbf{V}\|_F^2+ \|\mathbf{U}\|_F^2),
\end{equation}
\noindent
with a linear operator $\mathcal{A}$ that maps $\bb{U}\bb{V}^T$ to $\mathcal{F}_{\Omega}\{ \Sigma_{i=1}^k \mathbf{s}_j^i \odot \mathbf{x}^i \}$. As shown in \cite{Holme_Rosenzweig_Ong_Wilke_Lustig_Uecker_2019}, the above formulation can be further relaxed to a nuclear norm minimization problem in terms of a matrix $\mathbf{Z}$,

\begin{equation}\label{enlive:nl}
    \arg\min \limits_{\mathbf{Z}} \|\mathbf{y}-\mathcal{A}(\mathbf{Z})\|_2^2 + 2\alpha \|\mathbf{Z}\|_*.
\end{equation}

\noindent
The authors claim that this connection provides an explanation for why ENLIVE favors low-rank solutions. However, it is worth noting that due to the computational cost associated with the nuclear norm minimization, the convex formulation \eqref{enlive:nl} is not used for practical computations. Instead, the non-convex formulation \eqref{ENLIVE} is employed to solve the auto-calibrated reconstruction problem. Notably, neither of the formulations in \eqref{JSENSE} nor \eqref{NLIVE} explicitly exploits sparsity in the object being imaged.

\subsubsection{The relationship between $L$, $C$ and $n$}

The primary objective of parallel imaging (PI) techniques is to reduce imaging time, which is achieved by acquiring fewer data points per coil and taking multiple coil measurements at the same time. From compressive sensing theory \cite{candes_romberg_tao_2006,recht_fazel_parrilo_2010}, it is well known that the incoherence in the sensing matrix is crucial for the successful solution of the problems in \eqref{sparseMRI} and \eqref{sparseSENSE}. And it is well accepted in the literature on pMRI that increasing the number of channels $C$ improves the incoherence of the sensing matrix $\bb{E}$. 
\\

\noindent
However, to understand the role of $C$ and $L$ versus the sparsity level $n$, a thorough analysis of the sensitivity encoding, partial Fourier sampling, and the sparsifying transformation in the overall measurement ensemble is required. As with linear inversion for solving \eqref{pMRI}, one would expect a threshold for the maximum reduction factor and a minimum number of coils for the compressive sensing approaches (e.g. $\ell_0$, $\ell_1$ or mixed norm minimization) to be successful. The problem becomes even more challenging, when the sensitivity encoding is unknown as in the auto-calibration framework.
\\

\begin{figure}[ht!]
    \centering
    \begin{minipage}{0.45\linewidth}
        \centering
        \includegraphics[width=\linewidth]{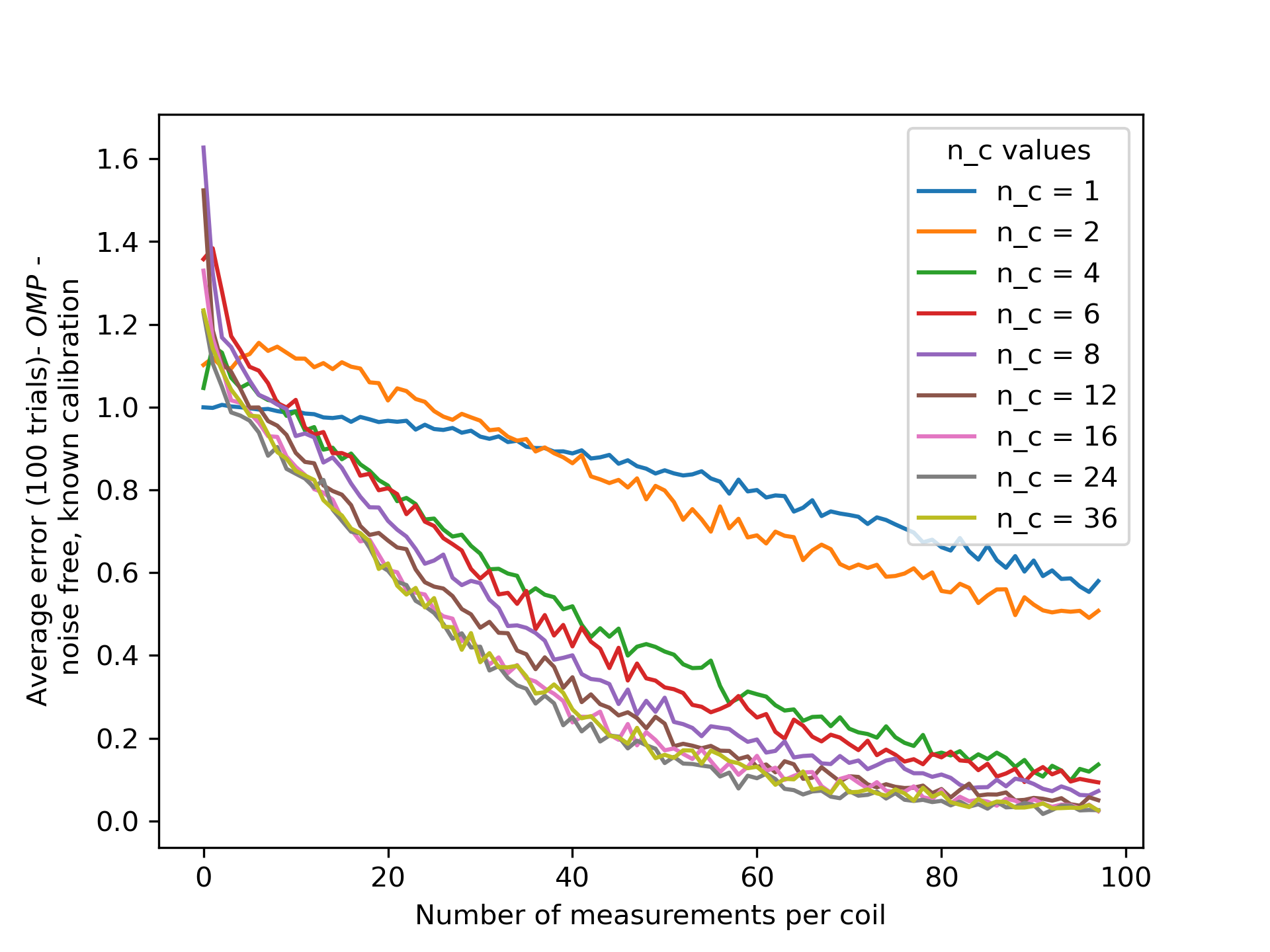}
        \caption*{(a) Noise-free case}
    \end{minipage} \hfill
    \begin{minipage}{0.45\linewidth}
        \centering
        \includegraphics[width=\linewidth]{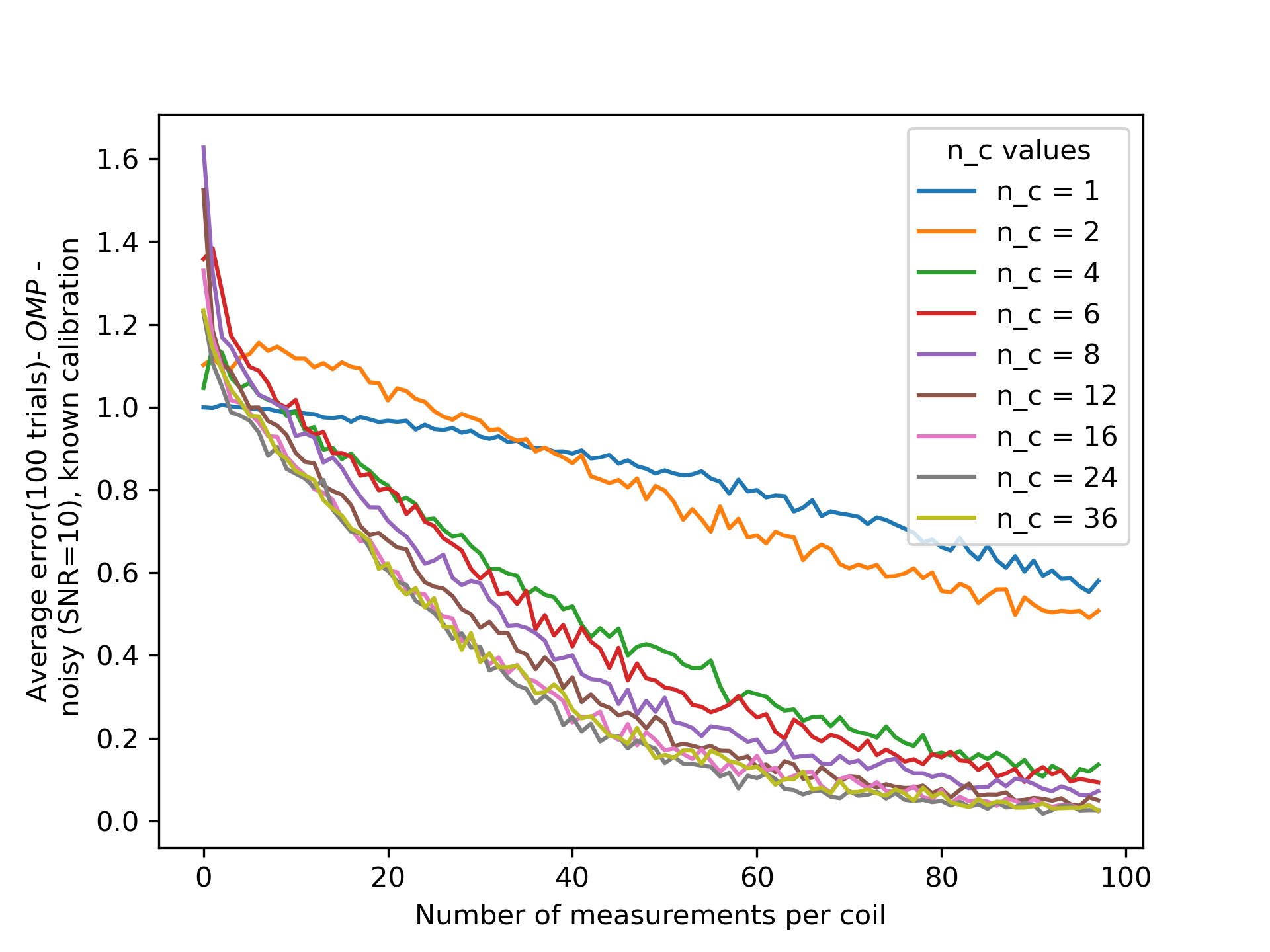}
        \caption*{(b) With independent Gaussian noise}
    \end{minipage}
    \caption{Two-dimensional pMRI simulations were performed, assuming truly transform n-sparse signals with randomly distributed non-zero components (N = 529, n = 32). Coils with sensitivities calculated using the Biot-Savart law were simulated for various numbers of coil elements ($n_c = 2, 4, 6, 8, 12, 16, 24,36$). In the case of a single coil, a constant sensitivity was assumed. The multi-coil signal was generated by multiplying the signal with the coil sensitivities, with optional addition of Gaussian noise. The signal is assumed to be sparse after 2-dimensional DCT transforms. The average of the relative $\ell_2$ error was computed over 100 different random undersampling realizations in Fourier-space and random realizations of transform sparse signals, using varying numbers of measurements and coil elements.}
    \label{fig:Otazo}
\end{figure}

\noindent
In this context, the works by Otazo et al. \cite{otazo2009distributed,Otazo_Kim_Axel_Sodickson_2010} should be mentioned,where they numerically assess the sufficient conditions on the number of measurements per coil $L$ and the number of coils $C$ for solving the $\ell_0$ norm minimization problem in \eqref{l0} or the $\ell_1$ norm minimization problem in \eqref{sparseMRI} to recover an $n$-sparse signal. The authors conclude that in the known calibration setup, when the number of coils is large enough, the required number of measurements per coil is very close to the signal sparsity for both the noiseless and noisy cases. To illustrate the numerical findings, we reproduced the reconstruction results from \cite{otazo2009distributed} by solving the $\ell_0$ minimization problem \eqref{l0}, using the orthogonal matching pursuit (OMP) algorithm. The results can be found in Figure \ref{fig:Otazo}. 

\vspace{10pt}

\noindent
At first glance, it may seem surprising that when $C$ is sufficiently large, increasing the number of coils does not yield any additional benefit in reducing the number of measurements per coil or the reconstruction quality. The following work may give us theoretical insights into the observation. Consider a more general setup in solving the following distributed compressed sensing problem:
\begin{equation}
    \mathbf{y}_i = \mathbf{F}_{\Omega_i}(\bb{s}_i \odot \bb{x}), i = 1 \cdots C, |\Omega_i|=L,
\end{equation}
where $\bb{s}_i$ are known and $\mathbf{x}$ is n-sparse.
\\

\noindent
Note that a key difference between the above problem and the pMRI forward model is that the sensors use different random undersampling patterns. The authors in \cite{Duarte_Sarvotham_Baron_Wakin_Baraniuk} have theoretically shown that, in such a distributed compressed sensing model with a large number of sensors, the required number of samples per sensor approaches the sparsity level, as solved by a greedy OMP-type algorithm. Recall that compressed sensing with \( \ell_1 \)-norm minimization (i.e. solving equation \eqref{sparseMRI} with a single coil) typically requires $L$ to be three to five times the number of sparse coefficients $n$ (i.e., \( 3 \)–\( 5n \) where $\bb{x}$ is n-sparse) to achieve accurate reconstructions \cite{lustig_donoho_pauly_2007,Tsaig_Donoho_2006}. Hence, a large number of sensors provides a constant improvement in the number of measurements per sensor, up to the sparsity level of the signal $\bb{x}$.
\\

\noindent
Note that this condition does not apply to multi-coil acquisition in pMRI, since all coils share the same undersampling pattern, albeit with distinct coil sensitivity profiles. Hence, one might expect a reduction in incoherence and a deviation from the theoretical recovery bounds discussed in \cite{Duarte_Sarvotham_Baron_Wakin_Baraniuk}. Additionally, the measurement ensemble in the pMRI setting must adhere to additional physical constraints, which makes the analysis within the compressive sensing framework more complex. By combining the numerical experiments with the aforementioned theory for the known calibration case, one could reasonably expect that a sufficient condition for robust reconstruction would adopt a similar form to that of the best-case scenario for the more complex auto-calibration problem in \eqref{eq:1}.

\subsubsection{Further discussions} 
Essentially, auto-calibration methods depend on data redundancy. In the literature on pMRI, other sources of redundancy have been utilized for auto-calibration, each under different assumptions. Data driven auto-calibration methods, such as GRAPPA \cite{griswold_jakob_heidemann_nittka_jellus_wang_kiefer_haase_2002} and SPIRiT \cite{Lustig_Pauly_2010} can be viewed as interpolation methods by estimating linear relationships within the fully sampled k-space (e.g., kernel calibration) and enforcing that relationship to synthesize data values in place of unacquired data (e.g., data reconstruction). More recently, a data-driven method SAKE \cite{shin_larson_ohliger_elad_pauly_vigneron_lustig_2014} structures the multi-coil dataset into a new data matrix that is designed to have low-rankness property. Then it solves a matrix completion problem to fill in the missing data. Joint estimation techniques attempt to iteratively estimate both the coil sensitivities and image contents while imposing some smoothness constraints on the sensitivity profiles, some examples include \cite{Ying_Sheng_2007,Uecker_Hohage_Block_Frahm_2008}. It is worth noting that the effectiveness of these methods depends on whether the coil and signal adhere to the assumptions made. As a result, they may produce varying quality of reconstructions for different datasets.

\subsection{Our contributions}
While the aforementioned auto-calibration methods in Section \ref{previous_methods} undoubtedly offer utility and applicability by giving successful numerical results, the majority of them do not provide rigorous analysis or theoretical guarantees of recovery. Note that even in the known calibration case \eqref{sparseSENSE}, the incoherence of the sensing matrix has not been properly analyzed for the popular $\ell_1$ minimization scheme, particularly for measurement ensembles that are applicable to pMRI. In this work, we take a step forward to provide analysis for the more complicated auto-calibration problem in \eqref{eq:1} and to provide sufficient conditions on $L$ and $C$ for stable and robust recovery of an approximate transform sparse signal $\bb{x}$.
\\

\noindent
We propose to solve the auto-calibration problem in \eqref{eq:2} using convex optimization by leveraging the block sparsity structure across all channels simultaneously. Our approach efficiently solves a mixed norm minimization problem and provides robust and stable recovery with theoretical guarantees. Interestingly, our theoretical findings align with previous numerical observations from past literature, providing valuable explanatory insights. To provide additional evidence, we present numerical results using the proposed method to solve the auto-calibrated pMRI. We include simulations using the analytical phantom and simulated coil sensitivities, as well as reconstructions from real pMRI measurements.

   \section{Problem setup}\label{3}  
        We are concerned with the following auto-calibration problem:
\begin{equation}
    \textbf{y}_i = \textbf{F}_\Omega \textbf{D}_i\textbf{x} + \textbf{n}_i, i = 1 \cdots C ,
\end{equation}
\noindent
 where $\bb{y}_i \in \mathbb{C}^{L \times 1}$ are the measurements,  $\bb{F}_\Omega \in \mathbb{C}^{L \times N}$ is the partial Fourier matrix with its rows chosen from the N-point 2-dimensional discrete Fourier transformation (DFT) matrix which belongs to a subset $\Omega \subset [N]$ ,
$\mathbf{D}_i = \text{diag}(\bb{s}_i) \in \mathbb{C}^{N\times N}$ is the diagonal matrix of calibration parameters corresponding to the $i$-th channel/coil, 
$\textbf{x} \in \mathbb{C}^{N \times 1}$ is the signal of interest and $\textbf{n}_i \in \mathbb{C}^{N \times 1}$ is additive noise.
\\

\noindent
Our goal is to recover the common signal of interest $\mathbf{x}$ from the measurements $\mathbf{y}_i$ across all channels, while the calibration matrices $\mathbf{D}_i$ are unknown. As discussed in the previous Sections, such a task is impossible unless additional assumptions are made on the signal and/or the calibration parameters.
\\

\noindent
Hence, we adopt the following two assumptions on both the signal $\mathbf{x}$ and the calibration parameters $\textbf{D}_i$, as justified in the previous discussions. In particular, 
\begin{enumerate}
    \item  $\textbf{x}$ has a sparse representation: there is a known sparsifying transformation $\mathbf{\Psi} \in \mathbb{C}^{N\times N}$ such that $\mathbf{z} = \mathbf{\Psi} \mathbf{x}$ is n-sparse. The location of the $n$ active coefficient of $\bb{z}$ is unknown;
    \item The calibration matrices belong to a common and known subspace: there is a known tall matrix $\textbf{B}\in\mathbb{C}^{N \times k}$ with $k < N$ such that $\mathbf{D}_i=\text{diag}(\mathbf{B}\textbf{h}_i), \text{for all } i = 1\cdots C$. Here, $\mathbf{h_i} \in\mathbb{C}^{k} $ are unknown.
\end{enumerate} 
\noindent
Under these assumptions, the recovery problem is to find $\mathbf{h}_i , i = 1\cdots C$ and a sparse $ \mathbf{z}$ that satisfies 
\begin{equation}\label{eq:constraint}
    \mathbf{y}_i = \mathbf{F}_\Omega \text{diag}(\mathbf{B}\textbf{h}_i) \mathbf{\Psi}^*\mathbf{z} , i = 1 \cdots C.
\end{equation}
\\

\noindent
Next, we apply the idea of `lifting' to transform the bi-linear constraint in both $\bb{h}_i$ and $\bb{z}$ to a linear constraint with respect to their outer product $\bb{h}_i\bb{z}^T$.

    \subsection{Lifting}
        Using the lifting approach \cite{Candès_Eldar_Strohmer_Voroninski_2013,Ahmed_Recht_Romberg_2014,ling_strohmer_2015}, we can transform the biconvex constraint \eqref{eq:constraint} into the problem of recovering a block-sparse matrix from linear measurements.
\\
\noindent
In particular, denote $\bb{B}=\begin{bmatrix} \mathbf{b}_1^* \\ \vdots \\ \mathbf{b}_N^*\end{bmatrix} = \begin{bmatrix}
    \mathbf{c}_1 &  \cdots & \mathbf{c}_k
    \end{bmatrix} \in \mathbb{C}^{N\times k}$, $\mathbf{\Psi}^* =\begin{bmatrix} \mathbf{q}_1^T \\ \vdots \\ \mathbf{q}_N^T \end{bmatrix}=\begin{bmatrix}
    \boldsymbol{\psi}_1 &  \cdots & \boldsymbol{\psi}_N
    \end{bmatrix} \in \mathbb{C}^{N\times N}$.
\\
A little linear algebra yields that the measurements $\bb{y}_i \in \mathbb{C}^{L\times 1}$ from each coil obey:
\begin{equation}\label{lifting}
    \bb{y}_i  = \bb{F}_\Omega \begin{bmatrix}
        \bb{b}_1^* \bb{h}_i \bb{z}^T \bb{q}_1 \\
        \bb{b}_2^* \bb{h}_i \bb{z}^T \bb{q}_2
        \\
        \vdots
        \\
        \bb{b}_N^* \bb{h}_i \bb{z}^T \bb{q}_N
    \end{bmatrix}
    =
    \mathbf{F_\Omega} \underbrace{\begin{bmatrix} (\mathbf{\bar{q}}_1 \bigotimes \mathbf{b}_1)^* \\ \vdots \\  \mathbf{\bar{q}}_N \bigotimes \mathbf{b}_N)^* \end{bmatrix}}_{\mathbf{\Phi}} \text{vec}(\bb{h}_i \bb{z}^T) := \mathbf{F}_\Omega \mathbf{\Phi} \text{vec}(\bb{h}_i \bb{z}^T), i = 1 \cdots C,
\end{equation} where \begin{equation}\label{Phi}
    \boldsymbol{\Phi} \coloneq \begin{bmatrix} (\mathbf{\bar{q}}_1 \bigotimes \mathbf{b}_1)^* \\ \vdots \\  (\mathbf{\bar{q}}_N \bigotimes \mathbf{b}_N)^* \end{bmatrix}   \in \mathbb{C}^{N\times kN}.
\end{equation}
\\

\noindent
Defining $$\bb{Y} = \begin{bmatrix}
    \bb{y}_1 | \bb{y}_2 | \cdots | \bb{y}_C
\end{bmatrix},$$ $$\bb{H} = \begin{bmatrix}
    \bb{h}_1 | \bb{h}_2 | \cdots | \bb{h}_C
\end{bmatrix},$$ $$\bb{X} = \begin{bmatrix}
    \bb{X}_1 | \bb{X}_2 | \cdots | \bb{X}_C
    \end{bmatrix}
    =
    \begin{bmatrix}
        \text{vec}(\bb{h}_1 \bb{z}^T) |
        \text{vec}(\bb{h}_2 \bb{z}^T) | \cdots |
        \text{vec}(\bb{h}_C \bb{z}^T)
    \end{bmatrix} = \bb{z} \otimes \bb{H} \in \mathbb{C}^{kN \times C},$$ we can rewrite \eqref{lifting} in matrix form as follows:
\begin{equation}\label{lifted}
    \bb{Y} = \underbrace{\bb{F}_\Omega \mathbf{\Phi}}_{\bb{A}} \bb{X},
\end{equation}
where $\bb{Y} \in \mathbb{C}^{L\times C}$, $\mathbf{F}_\Omega \in \mathbb{C}^{L\times N}$, $ \mathbf{\Phi} \in \mathbb{C}^{N\times kN}$, $\bb{X} \in \mathbb{C}^{kN\times C}$. We refer to $\mathbf{A} = \mathbf{F}_\Omega\mathbf{\Phi}$ as the forward matrix or measurement ensemble. Now the problem is to find the matrix $\mathbf{X}$ satisfying the linear constraint described by \eqref{lifted}. To that end we will exploit the block sparsity of $\bb{X}$.
\\
\subsection{Block sparsity of $\mathbf{X}$}
Assume $\bb{z}$ is $n$-sparse (i.e., $\|\bb{z}\|_0 = n$) and that $\text{supp}(\bb{z})= S$.  For any $i$, $\text{vec}(\bb{h}_i \bb{z}^T)$ is block n-sparse with blocks of length $k$. More precisely, define the indices of the $j$-th block as  $T_j = \{(j-1)k+1,\cdots, jk\}, j \in [N]$, $\text{supp}(\bb{h}_i \bb{z}^T) = \cup_{j \in S} \{(j-1)k+1,\cdots, jk\}$.  Concerning the sparsity of $\bb{X}$, since its columns share the same indexes for supports, it consists of at most n consecutive blocks of rows that are non-zero. More precisely, $\text{supp}(\bb{X})= (\cup_{j \in S} T_j) \times [C]$.
\\

\noindent
Moreover, $\bb{X}$ is at most rank-k if $C$ is large. Methods that exploit the low-rankness include \cite{Akçakaya_Basha_Goddu_Goepfert_Kissinger_Tarokh_Manning_Nezafat_2011,trzasko_manduca_2011,Holme_Rosenzweig_Ong_Wilke_Lustig_Uecker_2019} and also the the ENLIVE formulation \eqref{ENLIVE} introduced in Section \ref{2.2}. Naturally, one could minimize over a linear combination $\|\mathbf{\tilde{X}}\|_{1,2} + \lambda \|\mathbf{\tilde{X}}\|_*, \lambda > 0$ in an attempt to simultaneously capture the low-rankness and sparse characteristics of $\mathbf{\tilde{X}}$. However, as demonstrated in \cite{Oymak_Jalali_Fazel_Eldar_Hassibi_2015}, this approach is no more effective than using a single norm in terms of the number of measurements required. Theoretically, only recovering for a single coil case by minimizing $\|\cdot\|_*$ requires $L = O(N+k)$. Also $l_{1,2}$ minimization is more computationally efficient than the nuclear norm minimization problem. Hence, we choose to only minimize with respect to $\|\cdot\|_{1,2}$, which (as we will prove) may already yield a solution with sparse structure and small nuclear norm. 
\subsection{Optimization problem}
Given $(k,N,C)$, for any $\mathbf{X}\in \mathbb{C}^{kN \times C}$, define 
\begin{equation}\label{l12}
    \|\mathbf{X}\|_{1,2} = \sum\limits_{j \in [N]} \| \mathbf{X}_{T_j \times [C]}\|_F,
\end{equation} which is the sum of the Frobenius norm of the $k$ by $C$ blocks of $\mathbf{X}$. 
\\

\noindent
Since it is well known that the $\ell_{1,2}$-norm promotes the block sparse structure in $\mathbf{X}$ \cite{eldar_kuppinger_bolcskei_2010,eldar_mishali_2009,flinth_2017}, we solve the following convex program
\begin{equation}\label{optimization_problem}
     \min  \|\mathbf{X}\|_{1,2}
     \quad \text{subject to }    \bb{Y} = \bb{F}_\Omega \mathbf{\Phi} \bb{X} .
\end{equation}
\\
If the measurement are noisy, i.e., $\bb{Y} = \bb{F}_\Omega \mathbf{\Phi} \bb{X} + \mathbf{N}$ with $\|\mathbf{N}\|_F \leq \sqrt{C}\sigma$, we solve the following convex program
\begin{equation}\label{main0}
     \min  \|\mathbf{X}\|_{1,2}
     \quad \text{subject to}    \|\bb{Y} - \bb{F}_\Omega \mathbf{\Phi} \bb{X} \|_F \leq \sqrt{C}\sigma .
\end{equation}
\noindent
We also denote the forward operator as $\mathcal{A}$ : $\mathbb{C}^{k \times N} \rightarrow \mathbb{C}^L$ , or via a matrix $\mathbf{A} \in \mathbb{C}^{L \times kN}$,  such that  $\mathbf{Y} = \mathcal{A}(X) = \textbf{A} \mathbf{X}$. In our problem, $\bb{A} =\mathbf{F}_\Omega\mathbf{\Phi} $.
\\
\\
Once we get a solution $\hat{\bb{X}}$ to the $\ell_{1,2}$ minimization problem \eqref{optimization_problem}\eqref{main0}, we calculate the average of the best rank-one estimations of its columns to retrieve the common underlying signal. Recall the columns of $\hat{\bb{X}} = \begin{bmatrix}
    \hat{\bb{x}}_1 & \cdots \hat{\bb{x}}_C
\end{bmatrix} $ and the truth $\bb{X}^0 = \begin{bmatrix}
    \bb{X}^0_{1} & \cdots \bb{X}^0_{C}
\end{bmatrix} =  
\begin{bmatrix}
\text{vec}(\bb{h}_1 \bb{z}^T) |
        \text{vec}(\bb{h}_2 \bb{z}^T) | \cdots |
        \text{vec}(\bb{h}_C \bb{z}^T)
  \end{bmatrix}$. 
We use the following Lemma from \cite{ling_strohmer_2015} to compute an average of the estimates from $\hat{\mathbf{X}}$ and to ensure that it is close to $\bb{z}$. 

\begin{lemma}
    For any solution $\bb{X} = \begin{bmatrix}
        \bb{x}_1 | \cdots |\bb{x}_C
    \end{bmatrix}$ and the lifted matrix $\bb{X}_0 = \begin{bmatrix}
        \text{vec}(\bb{h}_1 \bb{z}_0^T) | \text{vec}(\bb{h}_2 \bb{z}_0^T) | \cdots | \text{vec}(\bb{h}_C \bb{z}_0^T)  )
    \end{bmatrix}$, 
    let $\hat{\sigma}\hat{\bb{u}}\hat{\bb{v}}^T$ be the best rank-one Frobenius norm approximation of  $\sum\limits_i \bb{x}_i/C$ reshaped into a matrix $\hat{\mathbf{X}}$ $\in \mathbb{C}^{k\times N}$ in the column-wise order. If $\|\bb{X} - \bb{X}_0\|_F = \sqrt{C}\epsilon $, then there exists a scalar $\alpha$ and a constant $C_0$ such that,
    \begin{equation*}
        \|\bb{z}_0-\alpha^{-1}\hat{\bb{v}}\|\leq C_0\min(\epsilon/\|\bb{z}_0\|,\|\bb{z}_0\|) .
    \end{equation*}
\label{lemma1}
\end{lemma}

\noindent 
\\
The Lemma is a direct consequence of combining Corollary 3.3 from \cite{ling_strohmer_2015} and using norm equivalency. Also notice if for any $\alpha_0 \neq 0$, $\alpha_0\bb{u}_0$ and $1/\alpha_0\bb{v}_0$ is also a pair of solutions. However, the scalar ambiguity is not important for many applications, including pMRI. The reconstructed signal is retrieved as $\hat{\bb{x}} = \boldsymbol{\Psi}^* \hat{\bb{v}}$ .
\\

\noindent
In the aspect of solving the optimization problem \eqref{optimization_problem}\eqref{main0}, the following unconstrained form is often used,
\begin{equation}\label{unconstraint}
    \min_{\bb{X}} \frac{1}{2}\|\bb{Y}-\bb{A}\bb{X}\|_{2}^2+\lambda\|\bb{X}\|_{1,2} .
\end{equation}
A comprehensive review of the algorithms to solve problems of the form~\eqref{unconstraint} is not a focus of this paper. We refer to a few common approaches such as the forward-backward splitting (FBS) \cite{Combettes_Patrick_Wajs_Valerie_2005}, ADMM \cite{wang_yang_yin_zhang_2008} and FISTA \cite{Beck_Teboulle_2009} which have found numerous applications in the context of pMRI reconstruction. In the simulations part of this paper, we use a FISTA-type iterative method to solve the unconstrained problem. Detailed iterative updates are provided in the Appendix \ref{Appendix A}. 

    \section{Main result}\label{4}
        Given any integer $k$ and $N$, we define $T_j \subset [kN]$ by $T_j = \{(j-1)k+1,\cdots,jk\}, j \in [N]$. Given an index set $S=\{s_1,s_2,\cdots\} \subset [N]$, define $S_j = T_{s_j}$ for $j \in |S|$ and $T := \cup_{j \in |S|} S_j$. With slight abuse of notation, we will not distinguish $\mathcal{P}_S$, $\mathcal{P}_T$ and $\mathcal{P}_{T\times [C]}$ when it is clear from the context. For example, for $\bb{x}\in \mathbb{C}^{N}$, $\bb{x}\in \mathbb{C}^{kN}$ and $\bb{X}\in \mathbb{C}^{kN \times C}$, $\mathcal{P}_S$ means the projection onto $S$, $T$ and $T\times [C]$ respectively. We will also sometimes write the projections $\mathcal{P}_S $ as subscripts, for example, $\bb{x}_S$ or $\bb{X}_S$, when it is convenient. 
\\

\noindent
For any $\bb{X}\in \mathbb{C}^{kN \times C}$, define $\text{sgn}(\bb{X}) \in \mathbb{C}^{kN \times C}$ as a matrix formed by normalizing each block of $\bb{X}$ with length $k\times C$. To be precise, 
\begin{equation}
    \text{sgn}(\mathbf{X}) =
    \begin{cases}
        \frac{\mathbf{X}_{T_j \times [C]}}{\|\mathbf{X}_{T_j \times [C]}\|_F}, \text{if} \|\mathbf{X}_{T_j \times [C]}\|_F \neq 0
        \\
        \bb{0}, \text{otherwise}
    \end{cases}.
\end{equation}
\noindent
It will come in handy in the proofs to decompose a linear map to any given subset of its columns. To that end, for a matrix $\bb{A} \in \mathbb{C}^{L \times kN}$ and a given set $S$, we define the projections onto its blocks of columns induced by $S$ as
\begin{equation}
    \bb{A}_{S_j} (:,k)= 
    \begin{cases}
        \bb{A}(:,k), \text{if } k \in S_j
        \\
        \mathbf{0}, \text{otherwise}
    \end{cases} \text{ and }
    \bb{A}_S (:,k)= 
    \begin{cases}
        \bb{A}(:,k), \text{if } k \in T 
        \\
        \mathbf{0}, \text{otherwise}
    \end{cases}.
\end{equation}
\\

\noindent
\subsection{Assumptions and main theorem}\label{assump}
Before presenting the main theorem, we first summarize the assumptions on the forward model in \eqref{lifted} and on the signals for easy reference.

\begin{enumerate}
    \item[(A1).] The rows of $\bb{B}$ are chosen uniformly at random without replacement from the rows of a matrix $\bb{B}_0$ with $\bb{B}_0^*\bb{B}_0 = \bb{I}_k$. (Equivalently, the columns of $\bb{B}_0$ form a tight frame.)
    \label{random_B}
    \item[(A2).] $\mathbf{\Psi}$ is a known and fixed orthonormal basis with $\mathbf{\Psi}^*\mathbf{\Psi} = \mathbf{\Psi}\mathbf{\Psi}^* = \bb{I}_N$.  \label{random_Psi}
    \item[(A3).] The rows of $\bb{F}_\Omega$ are chosen uniformly at random without replacement from the DFT matrix $\bb{F}$ with $\bb{F}^*\bb{F} = \frac{N}{L} \bb{I}_N$.  \label{random_F}
    \item[(A4).] We analyze two probability models/distributions for the sensitivity parameters:
    \begin{itemize}
        \item \textbf{Sampled from some basis of $\mathbb{C}^k$}:  Each column $\bb{h}_l,l \in [C]$ is chosen independently and uniformly at  random with replacement from the columns of an orthonormal basis $\bb{W} \in \mathbb{C}^{k\times k }$. 
        \item \textbf{Complex spherical}:  Each spanning coefficients $\bb{h}_l \in \mathbb{C}^{k},l \in [C]$ is chosen independently and uniformly at random from the complex sphere $S_{\mathbb{C}}^{k-1}$
    \end{itemize}
\end{enumerate}

\noindent 
We briefly introduced some common choices of $\bb{B}$ and $\boldsymbol{\Psi}$ in Section \ref{2}, including the smoothness promoting polynomial or sinusoidal basis for $\bb{B}$ and DCT or wavelet transformation for $\boldsymbol{\Psi}$. Here, we further comment on our choice of assumptions (A1)-(A4). (A1) and (A3) are regarding the sensing matrix. (A3) is a standard assumption in compressive sensing theory concerning the random sampling scheme for k-space/Fourier space measurements. For (A1), $\bb{B}$ is chosen uniformly at random from all possible row permutations of a fixed matrix $\bb{B}_0$. The orthogonality of the columns of $\bb{B}$ is not affected by row permutations since $\bb{B}_i^*\bb{B}_j = \sigma(\bb{B}_i)^* \sigma(\bb{B}_j)$. With $\bb{B}$ being a matrix with randomly permuted rows, for any fixed $\bb{h}$, the coil sensitivity $\mathbf{s} = \bb{B}\bb{h}$ will be randomly permuted copies of some underlying $\bb{s}_0 = \bb{B}_0\bb{h}$.
The randomness in $\bb{B}$ helps in providing an average case analysis for each coil. Practically, if a coil
$\bb{s}=\bb{B}\bb{h}$ were deterministic, any zero entries in 
$\bb{s}$ (e.g. dead pixels or inactive responses) would result in a loss of information about 
the object of interest $\bb{z}$ at those corresponding locations. Assumption (A4) on $\bb{H}$ is critical for developing an average case analysis across coils to explain the improved results when more channels are included. 
\\

\noindent
Furthermore, we consider the underlying image $\bb{z} = \mathbf{\Psi} \bb{x}$ to be fixed and only impose some probability distribution on $\bb{h}_i$. That is to assume we are imaging a fixed object with varying coil sensitivities where the spanning coefficient of the coils follows a certain probability distribution. For the two choices on the probability models on $\bb{H}$, they all ensure incoherence between the coils and that the coil sensitivities have the same energy. Observe first that, the incoherence between coils $\bb{s}_i^*\bb{s}_j = (\bb{B}\bb{h}_i)^*(\bb{B}\bb{h}_j) = \bb{h}_i^* \bb{B}^*\bb{B}\bb{h}_j = \bb{h}_i^*\bb{h}_j$ is small for $k$ large, and the coils all have the same energy $\bb{s}_i^*\bb{s}_i = (\bb{B}\bb{h}_i)^*(\bb{B}\bb{h}_i) = \bb{h}_i^* \bb{B}^*\bb{B}\bb{h}_i = \bb{h}_i^*\bb{h}_i = 1$.
\\

\noindent
Before moving on to stating the theorem, let us discuss the reason for imposing such a probability model on $\bb{H}$. Consider the worst case scenario where all the $\bb{h}_i$ are identical, then the measurements $\bb{Y} = \bb{A}\bb{X}$ consist of identical columns. By including multiple channels, we are not gaining more information on $\bb{X}$ regardless of the number of coil used. The following proposition gives a formal statement of the worst case analysis. The proof follows exactly from Proposition~4.1 in~\cite{Eldar_Rauhut_2010}.
\begin{lemma}{(\text{Worst case analysis})}
    Suppose there exist a $\bb{X}_1\in\mathbb{C}^{kN \times 1}$ that the $\ell_{1,2}$ minimization fails to receiver from $\bb{Y}_1 = \bb{A}\bb{X}_1$. Then the $\ell_{1,2}$ minimization fails to recover $\bb{X}$ = $ \begin{bmatrix}
        \bb{X}_1 & \bb{X}_1 & \cdots \bb{X}_1
    \end{bmatrix}  \in\mathbb{C}^{kN \times C} $ from $[\bb{Y}_1,\cdots,\bb{Y}_1] = \bb{A}\bb{X}$.
\end{lemma}

\noindent
However, one would expect that including multiple channels will provide more information to aid the recovery of the shared underlying object of interest. In a pMRI setting, it has been shown that parallel imaging (pMRI) outperforms single-coil MRI in many cases, and improved incoherence in the sensing matrix is often used as an explanation for such improvement. To understand the role of including multiple channels, and from the worst case analysis, it suggests a need to impose some probability model on the sensing parameters in order to prove some average case results.
\\

\noindent
Next, we state our main theorem that provides a sufficient condition on robust and stable recovery.
\begin{theorem}
   For any arbitrary subset $S \in [N]$ with $\text{card}(S) = n$ and $\bb{X}_0 := \bb{z}_0 \otimes \bb{H}, \bb{z}_0 \in \mathbb{C}^{N \times 1}, \bb{H} \in \mathbb{C}^{k \times C} $, consider the linear map $\bb{A} \in \mathbb{C}^{L\times kN}$ as defined in \eqref{lifted} that satisfies the randomness assumptions (A1)-(A4) and the noisy measurement $\bb{Y} = \bb{A}\bb{X}_0 + \boldsymbol{\mathcal{N}}$ with $\|\boldsymbol{\mathcal{N}}\|_F \leq \sqrt{C}\sigma$. The solution $\bb{X}$ to the $\ell_{1,2}$ minimization problem \eqref{main0} satisfies
\begin{equation}\label{main}
    \|\bb{X}-\bb{X}_0\|_F \leq C_1 \|\mathcal{P}_{S^c}\bb{X}_0\|_{1,2} + (C_2 + C_3 \sqrt{s})\sqrt{C}\sigma
\end{equation}
with probability $1-4N^{-\alpha}$ if the following conditions are met,
\begin{enumerate}
    \item[(C1.)] the number of measurement satisfies
    \begin{equation*}
    L \geq C_\alpha kn {\mu^B_{max}}^2{\mu^\Psi_{max}}^2 \log(knN),
\end{equation*}

    \item[(C2.)]
        $\bb{h}_l$ satisfies the assumption (A4) and the number of coils satisfies $$C \geq C_\alpha\log(N)k.$$
\end{enumerate}

\end{theorem}
\noindent
Here, $\mu^B_{max} := \max\limits_{i,j}\sqrt{L}|\mathbf{B}_{i,j}|$, $\mu^\Psi_{max} := \max\limits_{i,j}\sqrt{N}|{\boldsymbol{\Psi}}_{i,j}|$, which measure the incoherence of the rows of $\bb{B}$ and of $\boldsymbol{\Psi}$, respectively. Note that this condition includes many useful examples--for instance, when $\boldsymbol{\Psi}$ is the unitary DFT matrix, the constant $\mu^\Psi_{max}$ achieves the optimal value of one \footnote{A relaxed version allowing the constants to grow with the logarithmic scale of $N$ is also possible.} The constant $C_\alpha$ grows linearly with respect to $\alpha$. 
\\

\noindent
A detailed bound for the constants in the theorem can be found in the proofs in Section \ref{6}. Next, we use Lemma \ref{lemma1} to get back the signal of interest from the overall lifted matrix $\mathbf{X}$, which also eliminates the $\sqrt{C}$ scaling in the error bound in \eqref{main}.
\\

\noindent
\textbf{Remark}:
The result provides a sufficient condition for stable and robust recovery in the presence of noise and sparse deficiency in the signal, requiring both the number of the measurements per coil and the total number of coils to be sufficiently large for the boundary case. It is important to emphasize that if the signal is truly sparse and there is no noise present, perfect reconstruction is achieved (i.e., the error is zero). Our theoretical bound exhibits the same form as the empirical findings in \cite{Otazo_Kim_Axel_Sodickson_2010} (even in their non-convex setup) and is validated by our numerical tests. Furthermore, it also implies that once $L$ and $C$ are both sufficiently large, increasing the number of coils does
not yield any additional benefit in reducing the number of measurements per coil or the reconstruction
quality. In contrast, as we will show in Lemma \ref{lemma8}, for small $C$, increasing the number of coils improves the probability of recovery. Furthermore, to understand the strength of the statement, note that if $S$ is the true support of an $n$-sparse $\bb{z}_0$, $\|\mathcal{P}_{S^c}\bb{X}_0\|_{1,2} =0$, then the number of required measurements $L$ is of the order of the sparsity level of $\bb{X}_0$, up to log factors (i.e., $kn ={{\|P_S(\bb{X}_0)\|}_0}$). For the constants in the bound of (C1), it will be meaningful when the magnitudes of the entries in $\mathbf{B}$ and $\boldsymbol{\Psi}$ do not vary too much. Optimally, when the columns of $\mathbf{B}$ and $\boldsymbol{\Psi}$ are chosen from the DFT matrices, we optimally have $\mu^B_{max} = 1$, $\mu^\Psi_{max} =1$.
\\

\noindent
Ideally, we would like to derive a proof for the setup when both $\bb{B}$ and $\mathbf{\Psi}$ are fixed, and the randomness comes from the support set $S$ of $\bb{z}$ alone. The result would be analogous to the conditioning of the random sub-matrices of a fixed sampling matrix, e.g. see Chapter 14 of \cite{foucart_rauhut_2015}. Yet, we encountered some (at least to us) unsolvable difficulties in finding a realistic bound for the incoherence measure between the sub-blocks of the sampling matrix $\bb{A}$. Nevertheless, we still provide a formal statement of the bound and some discussions in the Appendix  \ref{Appendix}.
\\

\noindent
The proof of the main theorem follows a well-established proof outline of compressive sensing \cite{foucart_rauhut_2015}. First, we refer to a sufficient condition from \cite{flinth_2017} to ensure a robust and stable recovery. Subsequently, we analyze the circumstances under which these conditions can be met. 
We show in Lemma \ref{lemma8} that under assumptions (A1)-(A4), the sufficient conditions for recovery are satisfied with a higher probability as the number of coils increases. This means that, in practice, parallel imaging is likely to perform better than single imaging recovery.
 
\subsection{A sufficient condition }
We state below a sufficient condition for solving the general the $\ell_{1,2}$ minimization problem for recovering block sparse signals.
\begin{lemma} (Uniqueness Results)\label{sufficient}
Consider any matrix $\mathbf{X}_0 \in \mathbb{C}^{kN \times C}$ and arbitrary set $S \subset [N]$. Let $\bb{A}$ be a linear operator from $\mathbb{C}^{kN}$ to $\mathbb{C}^L$, and the noisy measurements $\bb{Y}= \bb{A}\mathbf{X}_0 + \boldsymbol{\mathcal{N}}$ with $\|\boldsymbol{\mathcal{N}}\|_F \leq \sigma$.
\\
Suppose that 
\begin{equation}
    \|(\bb{A}_S)^* \bb{A}_S - \bb{I}_S\|_{2\rightarrow 2} \leq \delta, \max \limits_{j \in |S^c|}\| \bb{A}_S^* \bb{A}_{S^c_j} \|_{2\rightarrow 2} \leq \beta
\end{equation}
for some $\delta \in [0,1)$ and $\beta \geq 0$.
\\
Moreover, suppose that there exists a matrix $\mathbf{V} \in \mathbb{C}^{L 
\times C}$ such that the approximate dual certificate $\boldsymbol{\mathcal{Y}} := \bb{A}^*\bb{V}$ satisfies the following conditions
\begin{equation}
    \|\mathcal{P}_S\boldsymbol{\mathcal{Y}} - \text{sgn}(\bb{X}_0)\|_F \leq \eta, \max\limits_{j\in |S^c|}\|\mathcal{P}_{S^c_j}\boldsymbol{\mathcal{Y}}\|_F \leq \theta \text{ and } \|\bb{V}\|_F \leq \tau\sqrt{s} .
\end{equation}
Define $\rho = \theta +\frac{\eta\beta}{1-\delta}$, $\mu = \frac{\sqrt{1+\delta}}{1-\delta}$. If $\rho < 1$, then the solution $\bb{X}^*$ to the problem \eqref{optimization_problem} obeys
\begin{equation*}
    \|\bb{X}^*-\bb{X}_0\|_F \leq 2\mu\sigma +(1+\frac{\beta}{1-\delta})(\frac{1}{1-\rho})(2\| \mathcal{P}_{S^c} \bb{X}_0\|_{1,2}+ 2\eta\mu\delta + 2\tau\sigma \sqrt{s}) .
\end{equation*}
\end{lemma}
\vspace{1pt}
\noindent
The derivation of the lemma is a direct application of Lemma 7 of \cite{flinth_2017} which relies on the sub-differential of $\partial (\|\cdot\|_{(1,2)})$. Given $\bb{Z}\in\mathbb{C}^{kN \times C}$, the sub-differential $\partial_\bb{Z} (\|\cdot\|_{(1,2)})$ is the set 
\begin{equation*}
    \{\bb{V} \in \mathbb{C}^{kN \times C} | \bb{V}_{T_i} := \frac{\bb{Z}_{T_i \times [C]}}{\|\bb{Z}_{T_i \times [C]}\|_F}, \text{if } \|\bb{Z}_{T_i \times [C]}\|_F \neq 0, \|\bb{V}_{T_i}\|_F \leq 1, \text{otherwise}\}. 
\end{equation*}
In the following sections, we will prove that with the assumptions we have made, the parameters $\delta$ and $\beta$ will be small with high probability. We will also construct a dual certificate $\boldsymbol{\mathcal{Y}}$ with small $\eta$, $\theta$ and $\tau$. Most of the proof procedure follows well-established methods in the compressive sensing literature, such as those described in \cite{ling_strohmer_2015} and \cite{flinth_2017}.

\subsection{Conditioning of $\mathbf{A}_S$}
Next, we give the conditions on when our measurement matrix $\mathbf{A}$ satisfies the condition $\|\bb{A}_S^*\bb{A}_S - \bb{I}_S\|_{2\rightarrow 2} \leq \delta$. Recall our forward matrix $\mathbf{A}_S=\mathbf{\mathbf{F}}_\Omega \mathbf{\Phi}_S$, the challenges come from the fact that its rows are dependent. Instead of bounding the norm directly, we apply the triangular inequality and work on bounding the two terms in $\underbrace{\| \mathbf{\Phi}_S^*\mathbf{\bb{F}}_\Omega ^* \mathbf{\bb{F}}_\Omega\mathbf{\Phi}_S - \mathbf{\Phi}_S^*\mathbf{\Phi}_S\|_{2\rightarrow2}}_{(1)}+\underbrace{\|\mathbf{\Phi}_S^*\mathbf{\Phi}_S - \mathbf{I}_S\|_{2\rightarrow2}}_{(2)}$ separately. Indeed, even when the full Fourier space is sampled, i.e., $\mathbf{F}_\Omega = \mathbf{F}$, the first term (1) is zero and we would at the very least, require the condition (2) to be satisfied. 

\begin{lemma}\label{lemma4}
    For operator $\mathbf{\Phi}$ defined in \eqref{lifting} which follows the assumptions induced by $\boldsymbol{\Psi}$ and $\bb{B}$ from (A1)-(A4) \ref{assump},
    \begin{equation}
    \|\mathbf{\Phi}_S^*\mathbf{\Phi}_S - \bb{I}_S \|_{2\rightarrow2} \leq \eta
    \end{equation}
    with probability at least $1-N^{-\alpha}$ if
    \begin{equation}
        N \geq C_\alpha {\mu^B_{max}}^2{\mu^\Psi_{max}}^2 kn \log N/\eta^2,
    \end{equation}
    where $C_\alpha$ grows linearly with respect to $\alpha$. $\mu^B_{max} := \max\limits_{i,j}\sqrt{L}|\mathbf{B}_{i,j}|$, $\mu^\Psi_{max} := \max\limits_{i,j}\sqrt{N}|{\boldsymbol{\Psi}}_{i,j}|$.
\end{lemma}
The condition given is mild since in practice we can properly assume $N \gg kn$.
\\

\noindent
The probability is over the randomness in $\bb{B}$. The proof of the Lemma is derived from Lemma 4.3 of \cite{ling_strohmer_2015}, where the authors consider the case $\bb{B}$ to be fixed but the rows of $\mathbf{\Psi}_S$ are chosen uniformly at random with replacement from the DFT matrix. 
\\

\noindent
\textbf{Remark}: Due to the symmetry of $\mathbf{\Psi}$ and $\bb{B}$ in constructing the forward matrix $\mathbf{\Phi}$ \eqref{lifting}, the result of Lemma \ref{lemma4} also holds true when $\bb{B}$ is fixed but the rows of $\mathbf{\Psi}$ chosen uniformly from all possible row permutations of a fixed matrix.    
\\

\noindent
The next lemma is concerned with  the tail bound for $\| \mathbf{\Phi}_S^*\mathbf{\bb{F}}_\Omega ^* \mathbf{\bb{F}}_\Omega\mathbf{\Phi}_S - \mathbf{\Phi}_S^*\mathbf{\Phi}_S\|_{2\rightarrow2}$.

\begin{lemma}\label{lemma5}
Suppose that the rows of $\mathbf{F}_\Omega$ are chosen uniformly at random without replacement from $\mathbf{F}$, where $\mathbf{F^*F}=\mathbf{F}\mathbf{F}^*=\frac{N}{L}$. For a realization of $\mathbf{\Phi}_S$ with $\|\mathbf{\Phi}_S^*\mathbf{\Phi}_S - \mathbf{I}_S\|_{2\rightarrow2}\leq \eta_1$,
\[\|\mathbf{\Phi}_S^*\bb{F}_\Omega^*\bb{F}_\Omega \mathbf{\Phi}_S-\mathbf{\Phi}_S^*\mathbf{\Phi}_S\|_{2\rightarrow2}\leq \eta , \, \text{with probability at least} \,\, 1-\epsilon, \]
provided that 
\[
L \geq \max \Big\{ C_{4}^2 \frac{(1+\eta_1)}{\eta^2}, C_3 \frac{1}{\eta} \Big\} {\mu_{max}^B}^2 {\mu_{max}^\Psi}^2 kn \, \ln\Big(\frac{2C_2^2k^2n^2}{\epsilon}\Big).
\]
\end{lemma}

\noindent 
Combining Lemma \ref{lemma4} and Lemma \ref{lemma5}, we have a probabilistic bound for $\|\bb{A}_S^*\bb{A}_S -\bb{I}_S\|_{2\rightarrow 2}$. Note that both the bounds are useful only when $\bb{B}$ and $\mathbf{\Psi}$ are ``flat'', i.e., $\mu^{\bb{B}}_{max} \sim O(1)$ and $\mu^{\mathbf{\Psi}}_{max} \sim O(1)$.
\\

\noindent
The next lemma provides a bound for the parameter $\beta$ in $\max \limits_{j \in |S^c|}\| \bb{A}_S^* \bb{A}_{S^c_j} \|_{2\rightarrow 2} \leq \beta$ in Lemma \ref{sufficient}.
\begin{lemma}\label{lemma6}
Consider the forward matrix $\bb{A}$ defined in \eqref{lifting} and assume that $\|\bb{A}_S^*\bb{A}_S -\bb{I}_S\|_{2\rightarrow 2}\leq \eta$, 

\begin{equation}
    \max\limits_{j \in [|S^c|]}\|\bb{A}_S^*\bb{A}_{S^c_j}\|_{2\rightarrow 2} \leq \sqrt{\frac{1+\eta}{L}} \mu_{max}^\mathbf{\Psi}.
\end{equation}

\end{lemma}

\subsection{Dual Certificate}

Regarding the exact dual certificate, let $\bb{V} = \bb{A}_S(\bb{A}_S^*\bb{A}_S)^{-1}\text{sgn}(\bb{X}_0)$, we have $\boldsymbol{\mathcal{Y}} = \bb{A}^*\bb{A}_S(\bb{A}_S^*\bb{A}_S)^{-1}\text{sgn}(\mathbf{X}_0)$ is the exact dual certificate. The next lemma deals with bounding the parameter $\tau$ in $\|\bb{V}\|_F \leq \tau\sqrt{s} $ in Lemma \ref{sufficient}. The bound obtained in the following lemma does not depend on the randomness of $\bb{H}$ or $sgn(\bb{X}_0)$. 
\begin{lemma}\label{lemma7}
    Let $\bb{V} = \bb{A}_S(\bb{A}_S^*\bb{A}_S)^{-1}\text{sgn}(\bb{X}_0)$. Conditioned on $\|\bb{A}_S^*\bb{A}_S -\bb{I}_S\|_{2\rightarrow 2}\leq \eta$,
     \begin{equation}
            \begin{split}
                \|\bb{V} \|_F
                \leq
                \frac{\sqrt{1+\eta}}{1-\eta} \sqrt{n}.
            \end{split}
        \end{equation} 
\end{lemma}

\begin{lemma}\label{lemma8}
Suppose that either each column of $\bb{H}$ is chosen independently and uniformly at random with replacement from the columns of an orthonormal basis $\bb{W} \in \mathbb{C}^{k\times k }$ or that each column of $\bb{H}$ is chosen independently and uniformly at random from the complex sphere $S_{\mathbb{C}}^{k-1}$ (i.e., $h_l \in \mathbb{C}^k \text{ with } \|h_l\|=1$). Conditioned on $\|\bb{A}_S^*\bb{A}_S - \bb{I}\|_{2\rightarrow2}\leq \eta$, if  
    \begin{equation}
         L > 6 {\mu^\Psi_{max}}^2 \frac{1+\eta}{(1-\eta)^2 } n,
    \end{equation}
then $\max\limits_{j \in [|S^c|]} \|\mathcal{P}_{S^c_j}\boldsymbol{\mathcal{Y}}\|_F < 1/2$ with probability at least $1- 2k\exp(-\frac{3C}{28k})$.
\end{lemma}

\noindent
Equivalently, we can also state the result in Lemma \ref{lemma8} as $\max\limits_{j \in [|S^c|]} \|\mathcal{P}_{S^c_j}\boldsymbol{\mathcal{Y}}\|_F < 1/2$ with probability at least $1-N^{-\alpha}$ if
\begin{equation*}
    C \geq \Tilde{C}_\alpha \log(N) k,
\end{equation*}
where $\Tilde{C}_\alpha$ grows linearly with respect to $\alpha$. This is the same form we adopt in the statement of the main theorem.
\\

\noindent
\textbf{Remark}: Note that increasing the number of channels $C$ increases the probability of recovery. 
And when $C$ is sufficiently large, the failure probability for recovery is dominated by the other sufficient conditions required in Lemma \ref{sufficient} hence further increasing $C$ does not provide any additional benefits in terms of the probability of recovery. Furthermore, the lower bound on $L$ in Lemma \ref{lemma8} is automatically met due to the bound $L \sim O(kn\log(k^2n^2 N))$ induced by the condition $\|\bb{A}_S^*\bb{A}_S - \bb{I}\|_{2\rightarrow2}\leq \eta$ from Lemma \ref{lemma5}.

    \section{Simulations}\label{5}
        In the first part of the simulations, we demonstrate the sufficient conditions of $L$ and $C$ for exact recovery using simulated data with different underlying dimensions of $k$ and $n$. We calculate the empirical probability of success for different values of $k$ and $n$ and show that the transition curves for $L$ and $C$ are consistent with the analytically derived bounds. In the second part, we present reconstruction results using simulated data and real pMRI measurements.

\subsection{Empirical success rate for the proposed method }

\begin{figure}[!ht]
   \centering
\begin{subfigure}{.45\textwidth}
  \centering
  \includegraphics[width=\textwidth]{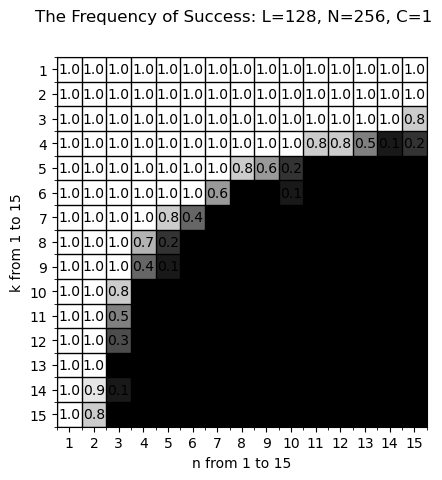}
  \caption{1 coil}
\end{subfigure}
\begin{subfigure}{.45\textwidth}
  \centering
  \includegraphics[width=\textwidth]{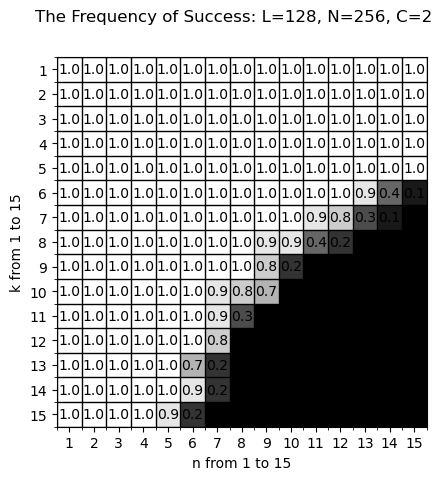}
  \caption{2 coil}
\end{subfigure}
\newline
 \begin{subfigure}{.45\textwidth}
  \centering
\includegraphics[width=\textwidth]{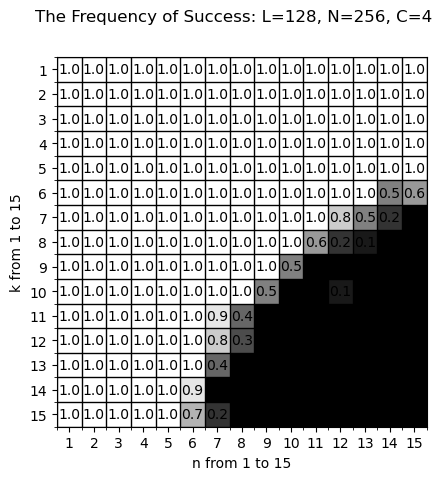}
  \caption{4 coil}
\end{subfigure}  
\begin{subfigure}{.45\textwidth}
  \centering
    \includegraphics[width=\textwidth]{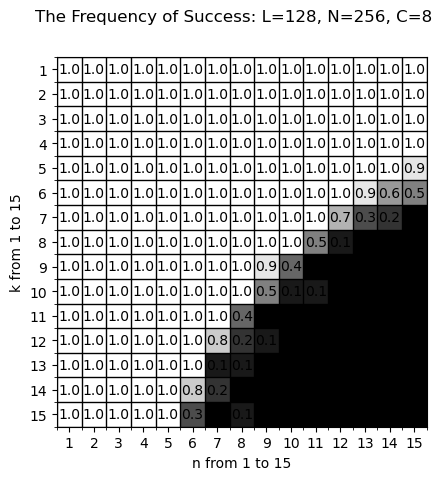}
  \caption{8 coil}
\end{subfigure}
\caption{Phase transition plot of performance by solving \eqref{main0} directly. The figures show empirical rate of success for a fixed sampling rate $L$ and different pairs of $(C,k,n)$. The numbers 1.0 means 100\% rate of success and 0.0 means 0\% rate of success. Observe that the transitional curve is improved for more $C$. However, when $C$ reaches a certain level the improvements in the empirical rate of success saturates.}  
\label{fig:2}
\end{figure}

\begin{figure}[!ht]
   \centering
\begin{subfigure}{.45\textwidth}
  \centering
  \includegraphics[width=\textwidth]{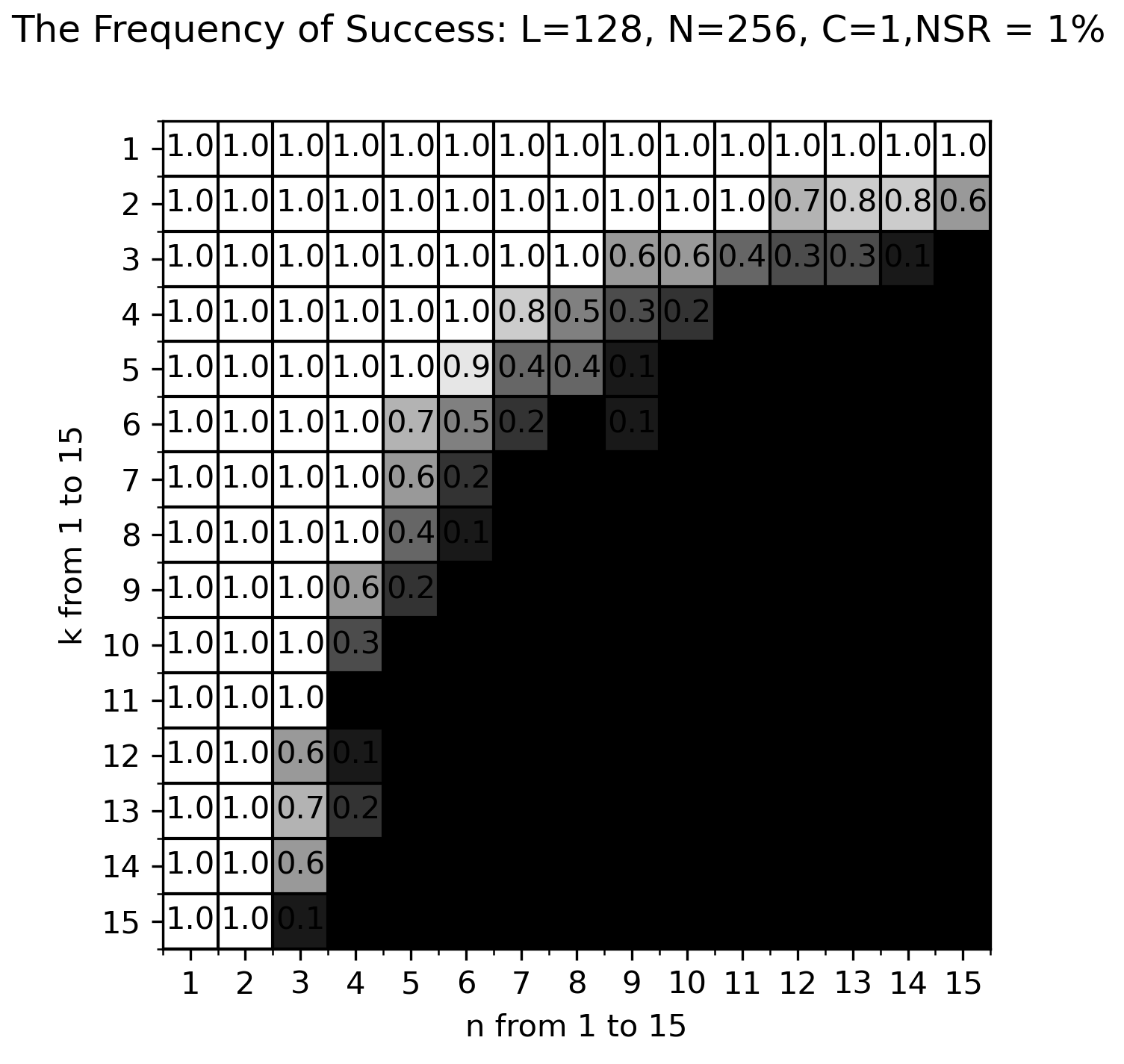}
  \caption{1 coil}
\end{subfigure}
\begin{subfigure}{.45\textwidth}
  \centering
  \includegraphics[width=\textwidth]{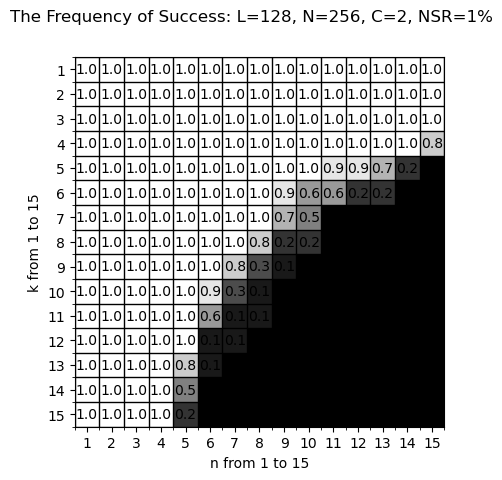}
  \caption{2 coil}
\end{subfigure}
\newline
 \begin{subfigure}{.45\textwidth}
  \centering
\includegraphics[width=\textwidth]{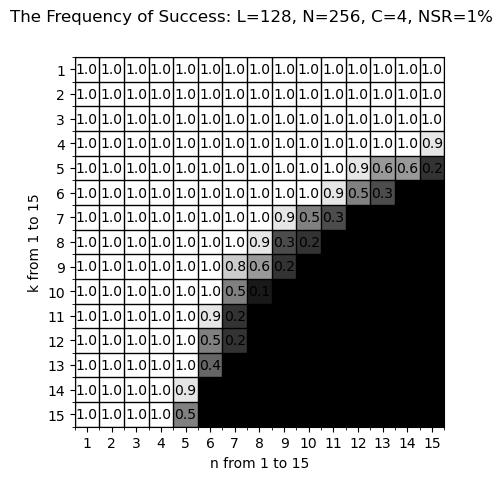}
  \caption{4 coil}
\end{subfigure}  
\begin{subfigure}{.45\textwidth}
  \centering
    \includegraphics[width=\textwidth]{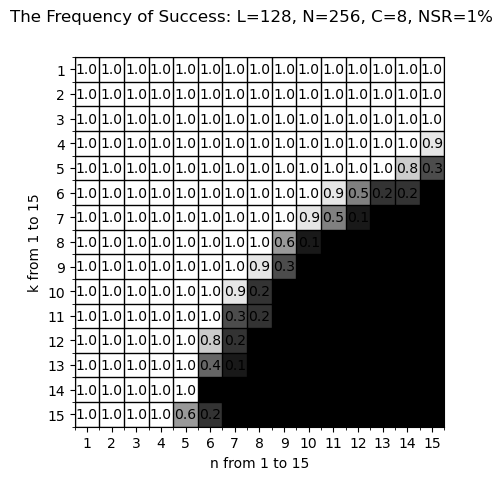}
  \caption{8 coil}
\end{subfigure}
\caption{Repeated the same experimental setup as in Fig:\ref{fig:2}, with 1\% of noise added.}  
\label{fig:3}
\end{figure}

For the first set of experiments, we follow the same setup as in Ling and Strohmer \cite{ling_strohmer_2015} and calculate the empirical probability of success of the proposed method. Fix the number of measurements $L = 128$ and the underlying signal dimension $N = 256$, and solve the optimization problem \eqref{main0} for varying $(k,n), k,n\in[15]$ and $C \in \{1,2,4\}$. For each $C$ and $(k,n)$, we generate $\bb{z} \in \mathbb{C}^N$ with random support with cardinality $n$ and the nonzero entries of $\bb{z}$ and $\bb{H} \in \mathbb{C}^{k \times C}$ drawn from $\mathcal{N}(0,1)$. Each time, we draw $L$ points uniformly from the Fourier space. We repeat the  experiment ten times, and each time new $\bb{H}$ and $\bb{z}$ are generated. 
\\

\noindent
We count an experiment as a success if the relative error between the solution of our optimization problem $\hat{\bb{X}}$ and the truth $\bb{X}_0$ is less than $0.01\%$, i.e., $\frac{\|\hat{\bb{X}}-\bb{X}_0\|_F}{\|\bb{X}_0\|_F} \leq 10^{-4} $. And the empirical success rate is a fraction defined as (total number of successes) / (total number of experiments). We use the cvxpy \cite{agrawal2018rewriting,diamond2016cvxpy} package of Python with the MOSEK \cite{mosek} solver to solve the mixed norm optimization problem \eqref{main0}.
\\

\noindent
Even though in the theorems we treated $\bb{B}$ as a random matrix with row permutations, in the experiments we use a fixed $\bb{B}$. We choose $\bb{B}$ as a fixed matrix consisting of the first $k$ columns of a random orthonormal matrix drawn from the $O(N)$ Haar distribution. It represents a stronger experiment than the theoretical treatment, as the theory only ensures recovery for most $\bb{B}$, while experimentally, we show the performance of a more relevant circumstance in practice. The results for choosing the DCT matrix as the sparsifying transformation (i.e., $\boldsymbol{\Psi}$) are shown in Fig:\ref{fig:2} and Fig:\ref{fig:3} for the noise-free and noisy cases, respectively.
\\

\begin{figure}[!ht]
   \centering
\begin{subfigure}{.45\textwidth}
  \centering
  \includegraphics[width=\linewidth]{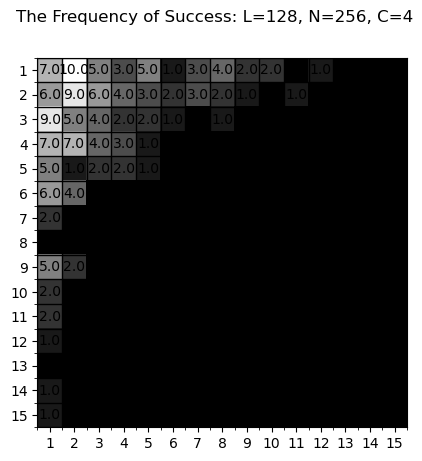}
  \caption{1 coil}
  \label{fig:dbsub1}
\end{subfigure}
\begin{subfigure}{.45\textwidth}
  \centering
  \includegraphics[width=\linewidth]{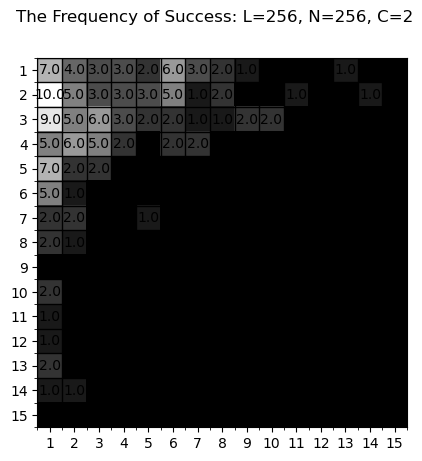}
  \caption{2 coil}
  \label{fig:dbsub2}
\end{subfigure}
\caption{Phase transition plot of performance by solving \eqref{main0} and using db4 wavelet transform as the sparsifying transformation matrix.}
\label{fig:4}
\end{figure}

\noindent
Next, we illustrate the barrier for reconstruction with the proposed framework when $\boldsymbol{\Psi}$ is a Wavelet transform. From the theory, a small mutual coherence, i.e., $\mu^\Psi_{max} \sim O(1)$, is required for the performance of the compressive sensing method. However, the Daubechies wavelet, which can be used to sparsify MRI images \cite{Chun_Adcock_Talavage_2014}, has high coherence. We followed the non-uniform sampling scheme in the Fourier measurements as suggested by \cite{Chun_Adcock_Talavage_2014,ADCOCK_HANSEN_POON_ROMAN_2017} as an attempt to improve the results. However, even when $L=N$ (when Fourier space is fully sampled), recovery is poor. We use the Daubechies 4 wavelet as an illustration, and the results are shown in Figure \ref{fig:4}. 

\vspace{10pt}
\noindent
Nevertheless, since our method is designed to work with approximately sparse signals, and given that we will demonstrate the effectiveness of using the Discrete Cosine Transform (DCT) as the sparsifying transformation, it is not critical if the wavelet does not perform as well in this context. 
\subsubsection{Improved results compared to $\ell_1$ minimization}
\begin{figure}[!h]
   \centering
\begin{subfigure}{.3\textwidth}
  \centering
\includegraphics[width=\linewidth]{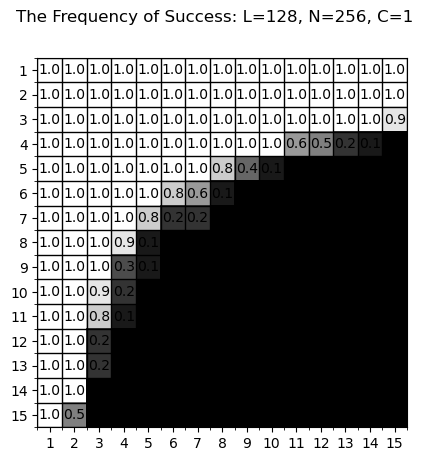}
  \caption{1 coil}
  \label{fig:l1sub1}
\end{subfigure}
\begin{subfigure}{.3\textwidth}
  \centering
  \includegraphics[width=\linewidth]{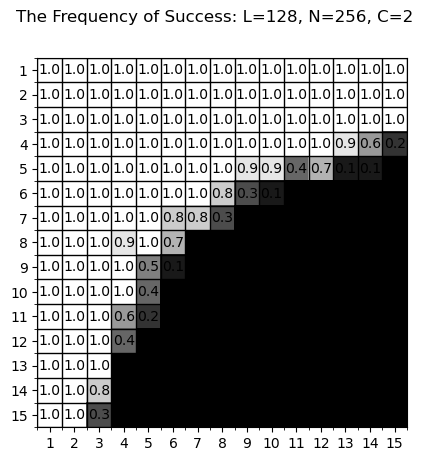}
  \caption{2 coil}
  \label{fig:l1sub2}
\end{subfigure}
\begin{subfigure}{.3\textwidth}
  \centering
  \includegraphics[width=\linewidth]{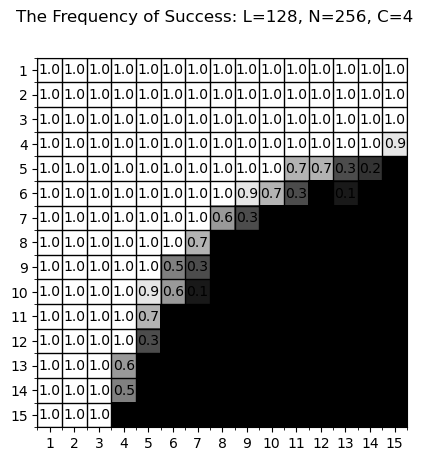}
  \caption{4 coil}
  \label{fig:l1sub4}
\end{subfigure}
\caption{Phase transition plot of performance by minimizing $\sum\limits_{i \in [C]}\|\bb{X}(:,i)\|_2$ instead. Observe that compared to the proposed method, the improvement in the success rate is not obvious for including larger $C$.} 
\label{fig:5}
\end{figure}

Our proposed method minimizes $\sum \limits_{i \in [N]} \|\textbf{X}_{S_i}\|_2$ to promote block sparsity within a lifted signal and between coils. We demonstrate its advantage over minimizing $\sum\limits_{i \in [C]}\|\bb{X}(:,i)\|_2$, which only promotes row sparsity between coils. We change the objective in \eqref{main0} to minimize $\sum\limits_{i \in [C]}\|\bb{X}(:,i)\|_2$ and use the same experimental setup. The results are shown in Figure \ref{fig:5}.

\subsubsection{Minimal L required for exact recovery is proportional to $kn$}
Next, we perform two sets of experiments to show that the minimal measurements $L$ required for exact recovery are proportional to the sparsity level $kn$. 
\\

\begin{figure}[!ht]
   \centering
\begin{subfigure}{.4\textwidth}
  \centering
  \includegraphics[width=0.9\linewidth]{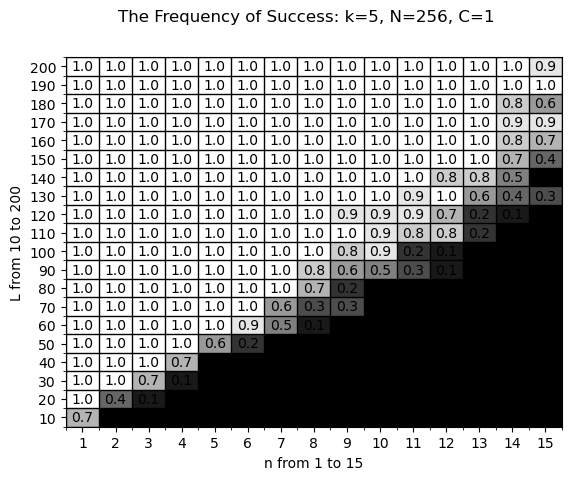}
  \caption{L prop n, 1 coil}
\end{subfigure}
\begin{subfigure}{.4\textwidth}
  \centering
  \includegraphics[width=0.9\linewidth]{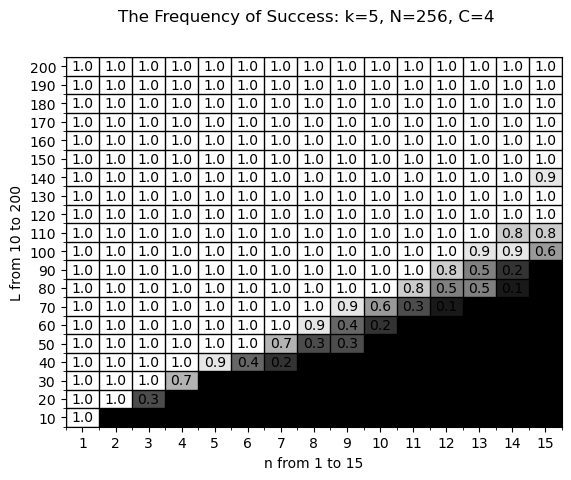}
  \caption{L prop n , 4 coil}
\end{subfigure}
\caption{The empirical rate of success for a fixed k and varying $(L,n,C)$.}   \label{fig:6}
\end{figure}

\begin{figure}[!h]
   \centering
\begin{subfigure}{.4\textwidth}
  \centering
  \includegraphics[width=0.9\linewidth]{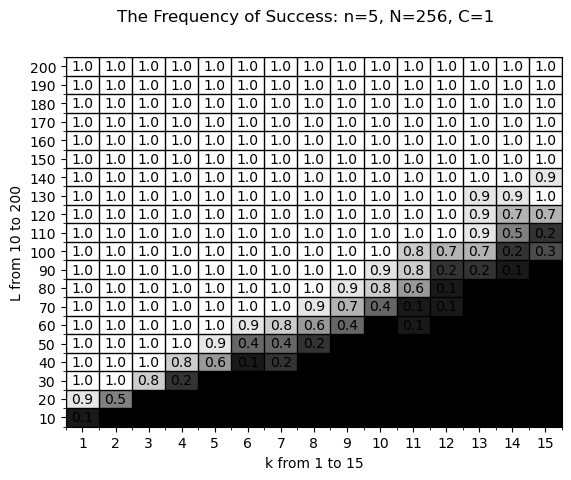}
  \caption{L prop k, 1 coil}
\end{subfigure}
\begin{subfigure}{.4\textwidth}
  \centering
  \includegraphics[width=0.9\linewidth]{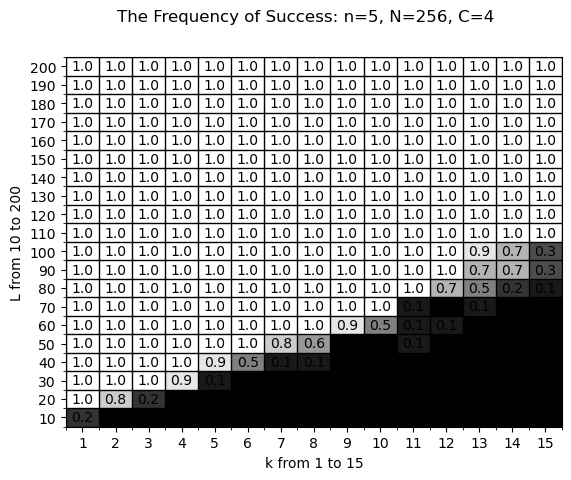}
  \caption{L prop k , 4 coil}
\end{subfigure}
\caption{The empirical rate of success for a fixed n and varying $(L,k,C)$. Observe that the transitional curves for all four figures shows a linear trend.} 
\label{fig:7}
\end{figure}

\noindent
\begin{itemize}
    \item For a fix $N=256$, $n=5$, we let $k$ ranges from 1 to 15 and $L$ varies from 10 to 200. We run the simulation for 10 times and calculate the empirical success rate.
    \item For a fix $N=256$, $k=5$, we let $n$ ranges from 1 to 15 and $L$ varies from 10 to 200.  We run the simulation for 10 times and calculate the empirical success rate.
\end{itemize} 
We also repeat the same experiments for different $C=1,4$. The results are shown in Figure \ref{fig:6} and \ref{fig:7}. 

\subsection{Reconstruction of the analytical phantom}
For the second part of the simulation, we show reconstruction results using the analytical phantom and simulated coils. We point out the major differences of the experimental setup from the last part:
\begin{itemize}
    \item A crucial difference in the simulations is that we do not manually enforce the object of interest (i.e. $\bb{x}$ or $\bb{z}$) to be sparse/transform sparse. Hence, it allows for sparsity deficiency.
    \item Different levels of noise are added to the measurements. 
    \item Instead of solving the constrained optimization problem \eqref{main0}, we solve the regularized problem \eqref{unconstraint} using an ISTA-type of updates as shown in Appendix \ref{Appendix A}.
\end{itemize}
\noindent
We use the analytical phantom \cite{Guerquin-Kern_Lejeune_Pruessmann_Unser_2012} of size $256$ by $256$ and randomly generated coil sensitivities following the Biot–Savart law. It is important to note that we do not manually enforce sparsity in our signal. For example, one may manually use a sparse estimation of the phantom $\bb{x}$, e.g., by crafting the signal by thresholding $\bb{z} =  \boldsymbol{\Psi}\bb{x} = 0 , \text{ if } |\boldsymbol{\Psi}\bb{x}| < \tau$ and set $\bb{x} = \boldsymbol{\Psi}^* \bb{z}$. The thresholding effectively reduces the dimensionality of the problem and would favor the reconstruction results. However, such manipulation is only possible in simulations when the true signal is known beforehand. Our approach is more flexible and allows for sparsity deficiency in the underlying signal.
\\

\noindent
We use different number of coils to show the efficiency of the proposed method. We plot reconstruction results using different levels of sub-sampling rate and Gaussian noise level (i.e., noise to signal ratio). We report the relative error (rerr) of the reconstructed image and the true phantom. The reconstruction results are shown in Figure \ref{fig:8}.
\newpage
\begin{figure}[hbt!]
\captionsetup{font=small}
 \begin{subfigure}{\textwidth}
 \centering
       \includegraphics[width=\linewidth]{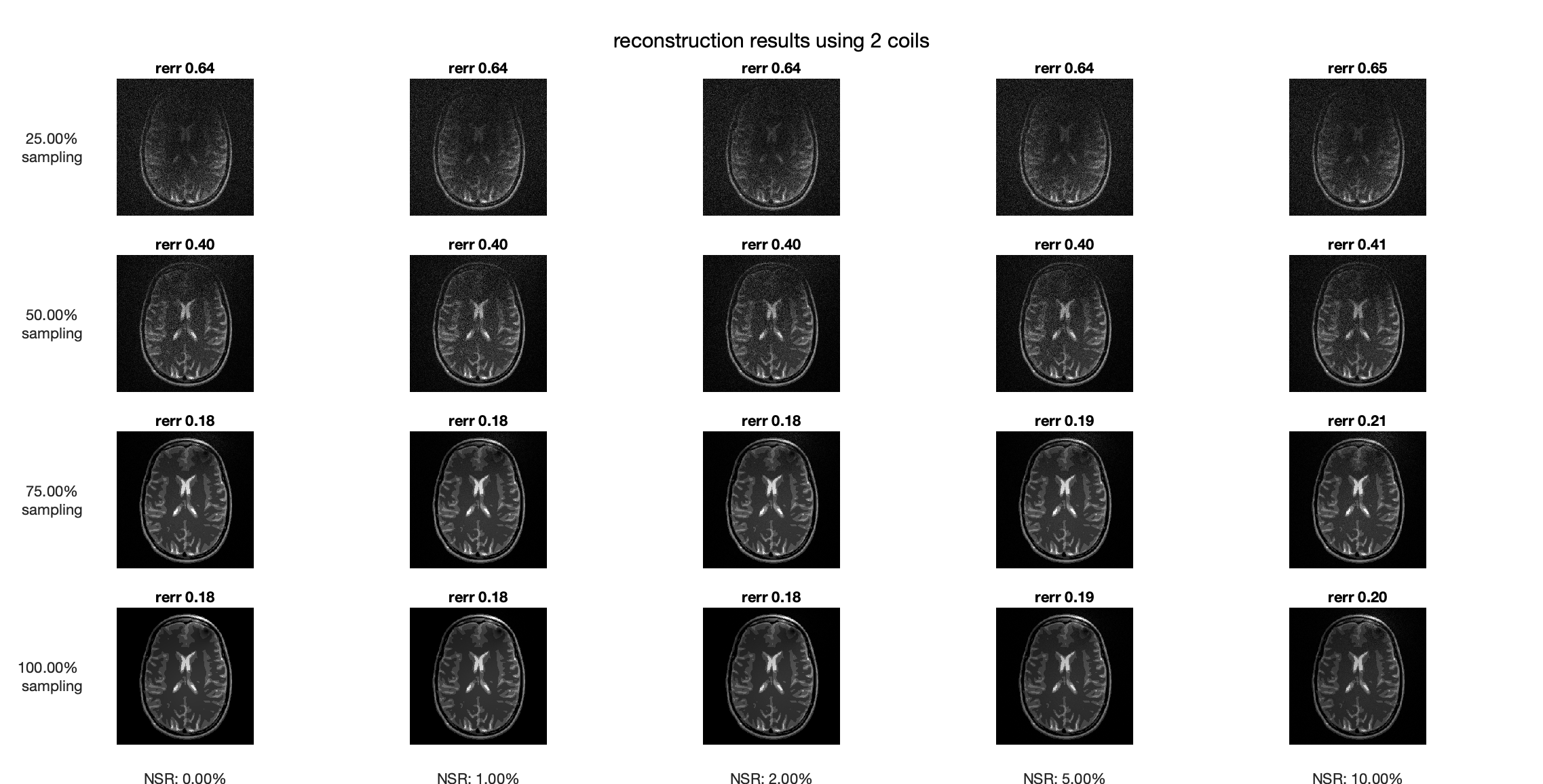}
  \caption{2 coils}
  \end{subfigure}
 \begin{subfigure}{\textwidth}
 \centering
       \includegraphics[width=\linewidth]{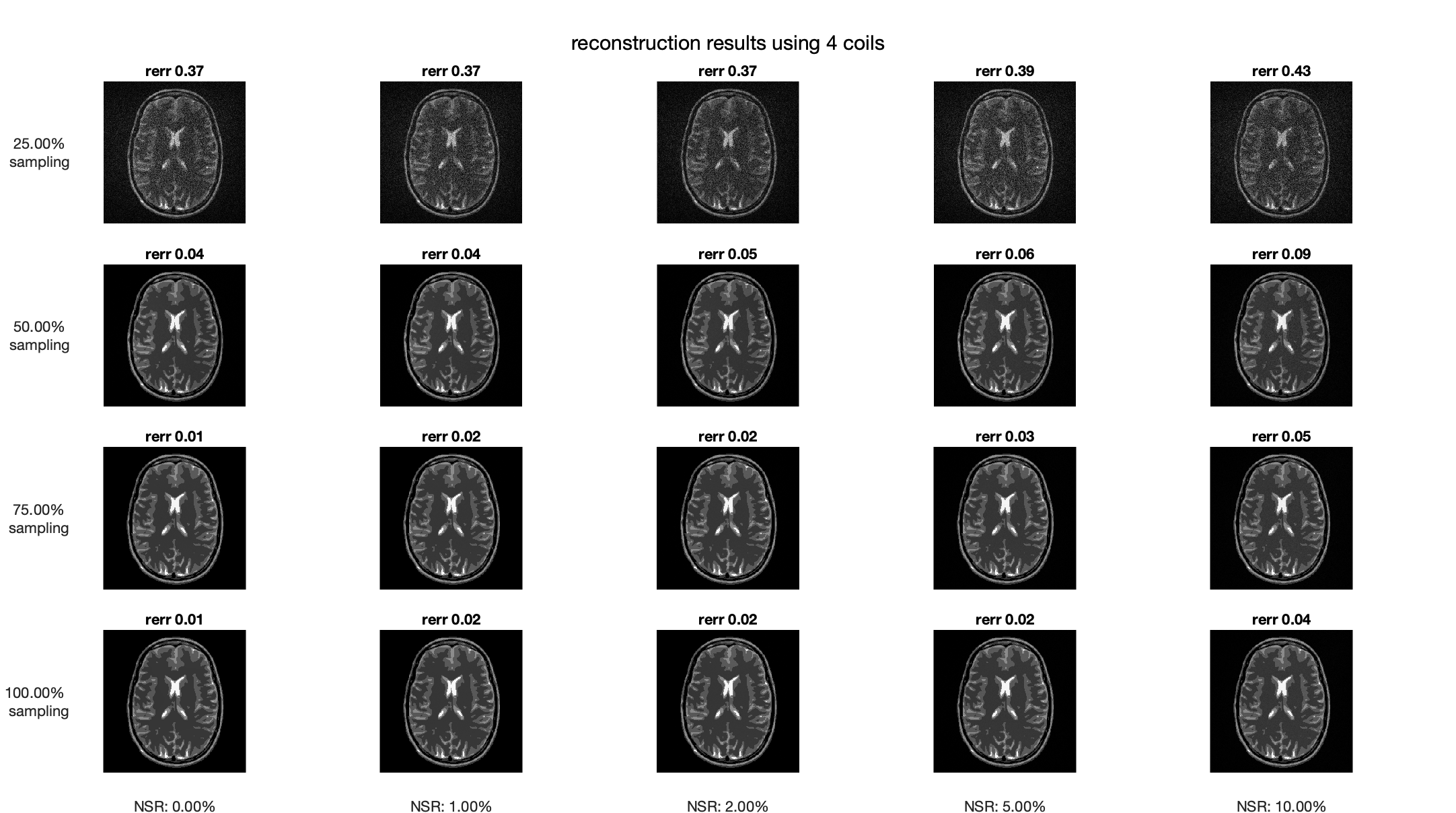}
  \caption{4 coils}
   \end{subfigure}
\end{figure}
\begin{figure}\ContinuedFloat
 \begin{subfigure}{\textwidth}
 \centering
       \includegraphics[width=\linewidth]{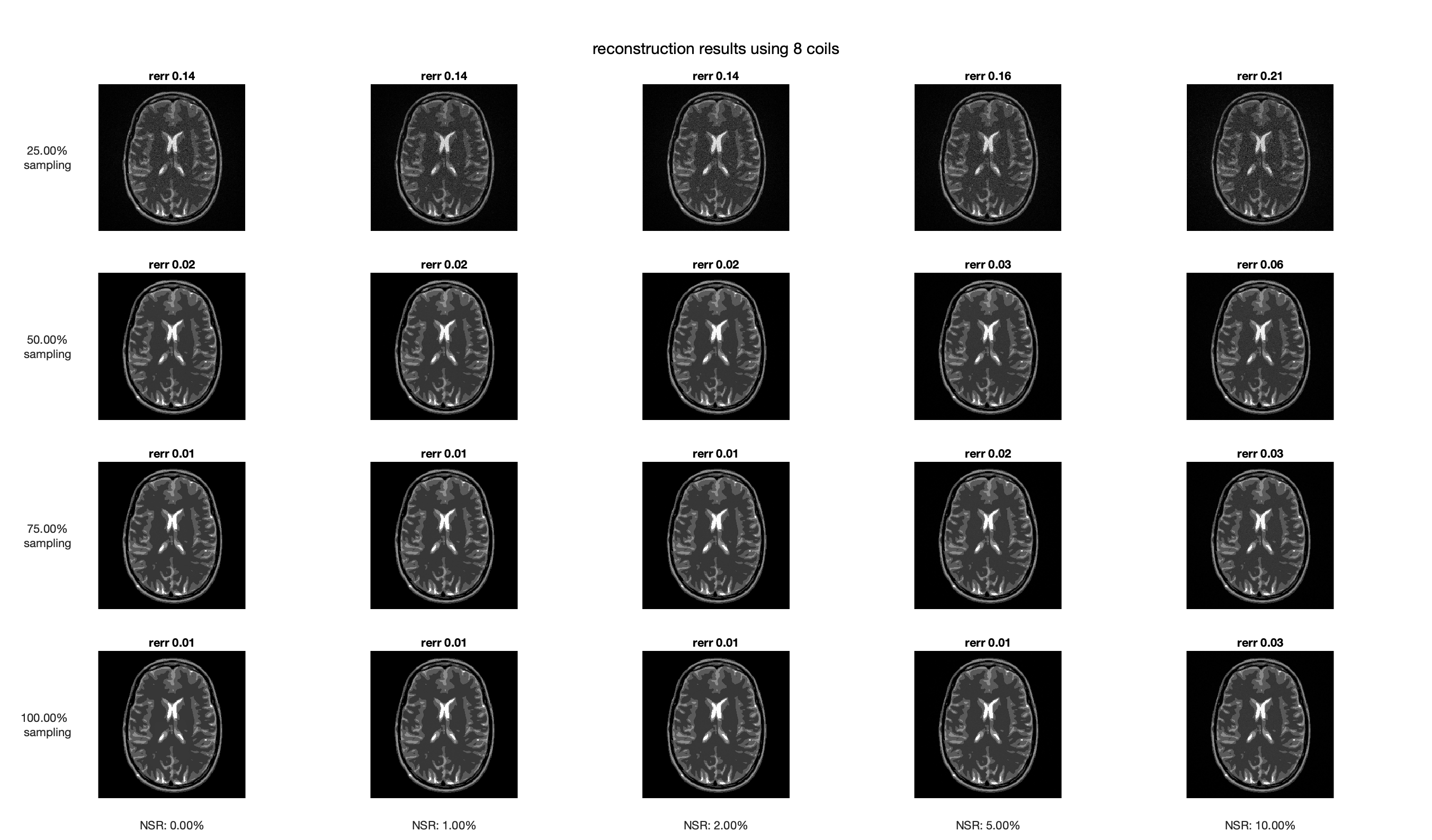}
  \caption{8 coils}
   \end{subfigure}
   
\caption{Reconstruction of the analytical phantom. The x-axis and y-axis represents varying noise-to-signal ratios and the sub-sampling rate respectively. Our observation indicates a gradual degradation in reconstruction quality with increased noise and increased reduction factor.}
\label{fig:8}
\end{figure}
\newpage
\noindent
To assess the performance of the proposed method with respect to the number of coils, we plot the average reconstruction error versus the number of measurements per coil, for different numbers of coils and levels of Gaussian noise. We adopted a similar setup as that in \cite{otazo2009distributed,Otazo_Kim_Axel_Sodickson_2010} where the multi-coil signal was generated by multiplying a small analytical phantom (N = 529) by the coil sensitivities (k=6) and optionally adding Gaussian noise. The average reconstruction error (relative $\ell_2$ distance) over 10 different random Fourier space undersampling realizations was computed for different numbers of measurements and coil elements. A DCT transform was employed to sparsify the image and no strict sparsity was enforced in the signal. The number of required measurements to achieve the noise floor reconstruction decreases with the number of coils for both the noise-free and the noisy case. These trends are consistent with previous numerical experiments for the known calibration pMRI reconstruction, e.g. \cite{otazo2009distributed}. The plots confirm that with a large number of coils, the number of measurements per coil sufficient to achieve a baseline level of reconstruction error is significantly smaller than the true dimensionality of the signal.

\begin{figure}[ht!]
\centering
\begin{subfigure}{0.45\textwidth}
    \centering
    \includegraphics[width=\linewidth]{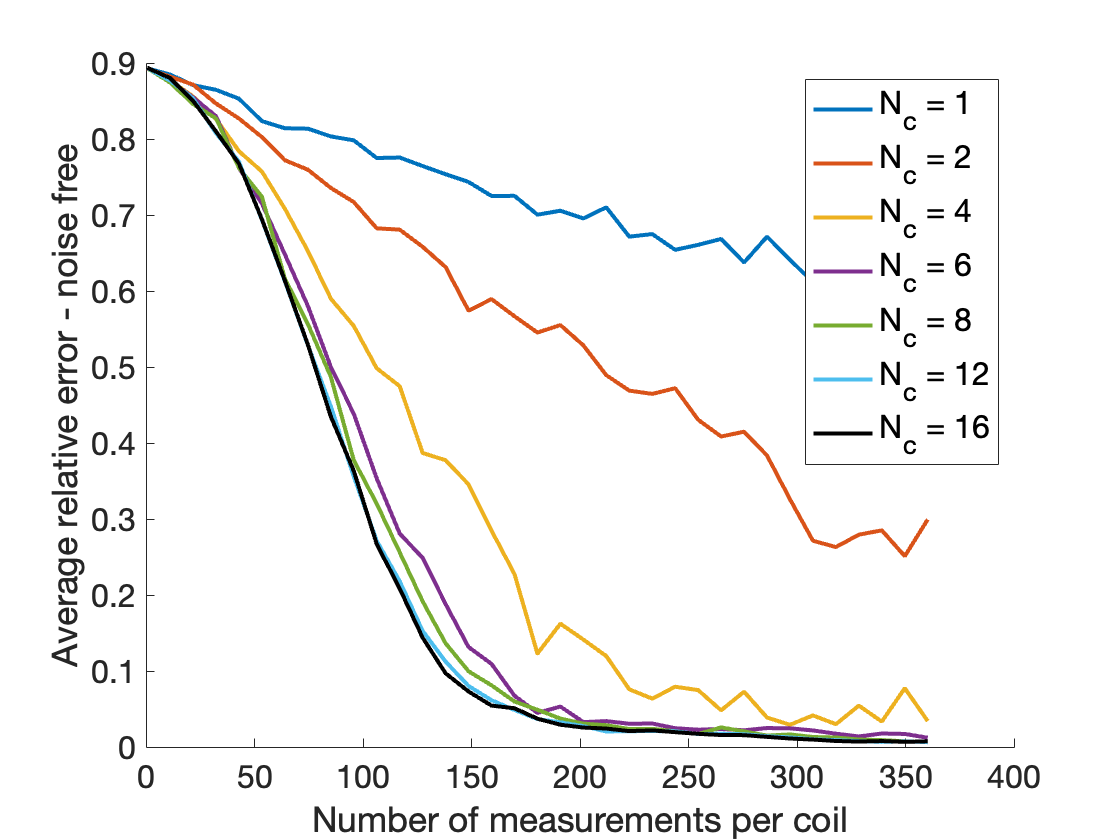}
    \caption{}
\end{subfigure}
\hspace{0.01\textwidth} 
\begin{subfigure}{0.45\textwidth}
    \centering
    \includegraphics[width=\linewidth]{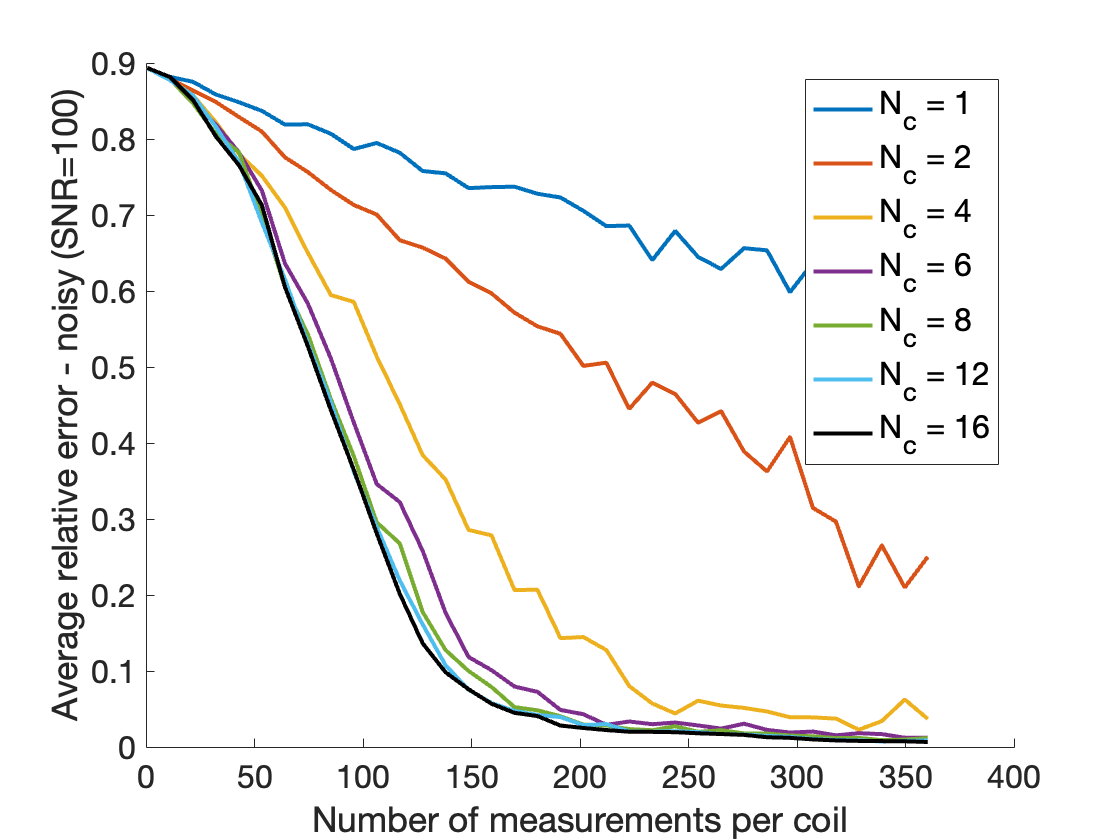}
    \caption{}
\end{subfigure}
\vspace{0.01\textwidth} 
\begin{subfigure}{0.45\textwidth}
    \centering
    \includegraphics[width=\linewidth]{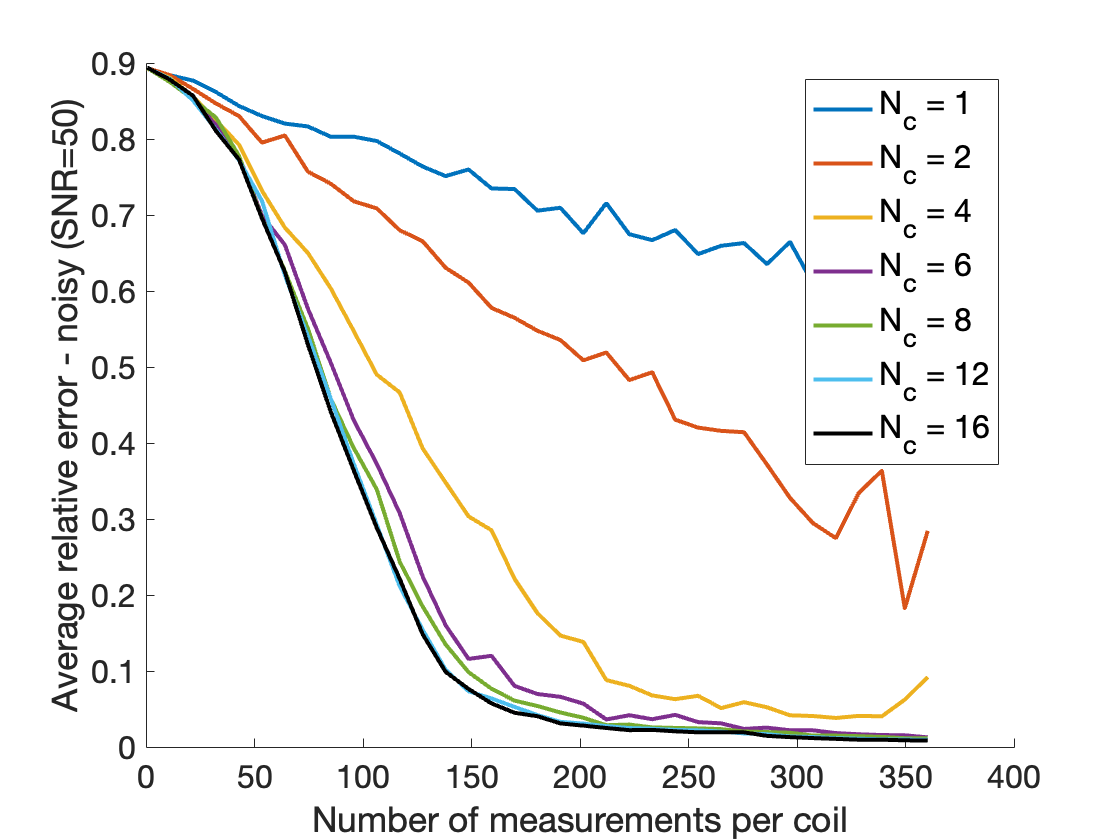}
    \caption{}
\end{subfigure}
\hspace{0.01\textwidth} 
\begin{subfigure}{0.45\textwidth}
    \centering
    \includegraphics[width=\linewidth]{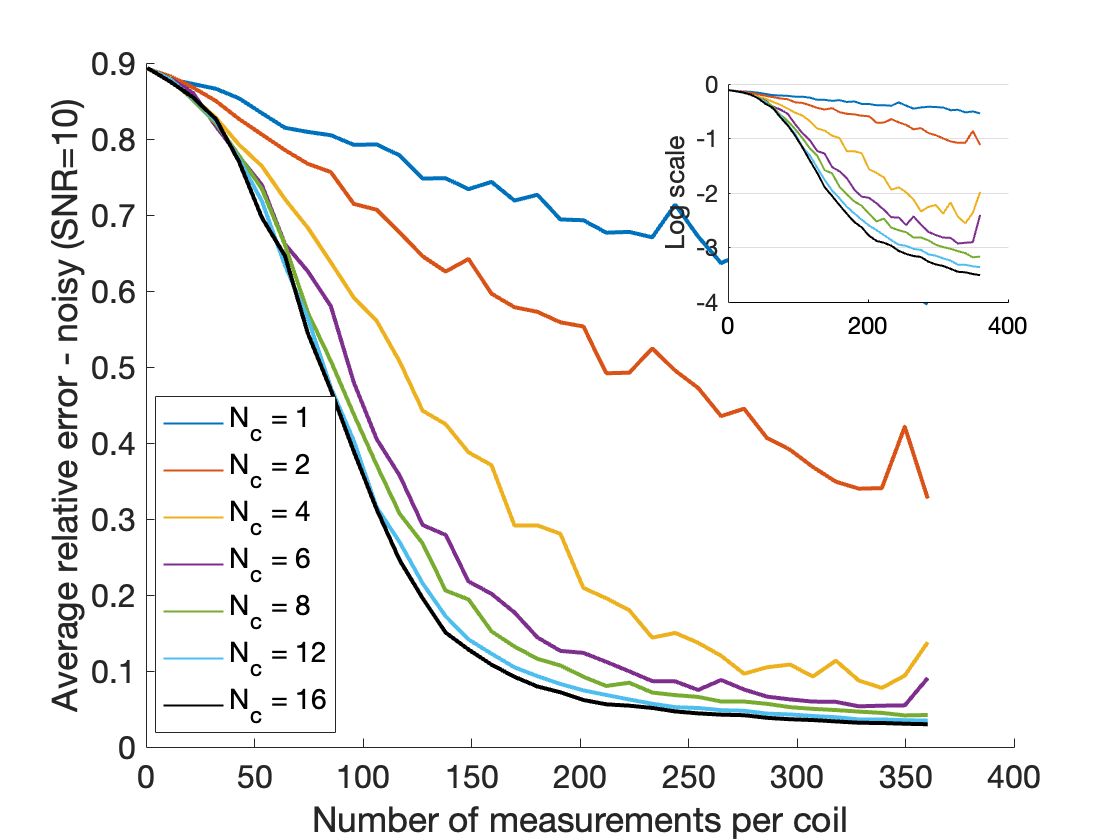}
    \caption{}
\end{subfigure}
\caption{The four subplots correspond to the noise-free case and SNR values of 100, 50, and 10. The relative $\ell_2$ error of the reconstruction is plotted against varying number of measurements per coil ($L$). For SNR=10, the log scale of the relative error is plotted in the inset of (d).}
\end{figure}

\subsection{Comparison with the ENLIVE method}
\begin{figure}[h!]
    \centering
    \includegraphics[width=0.75\linewidth]{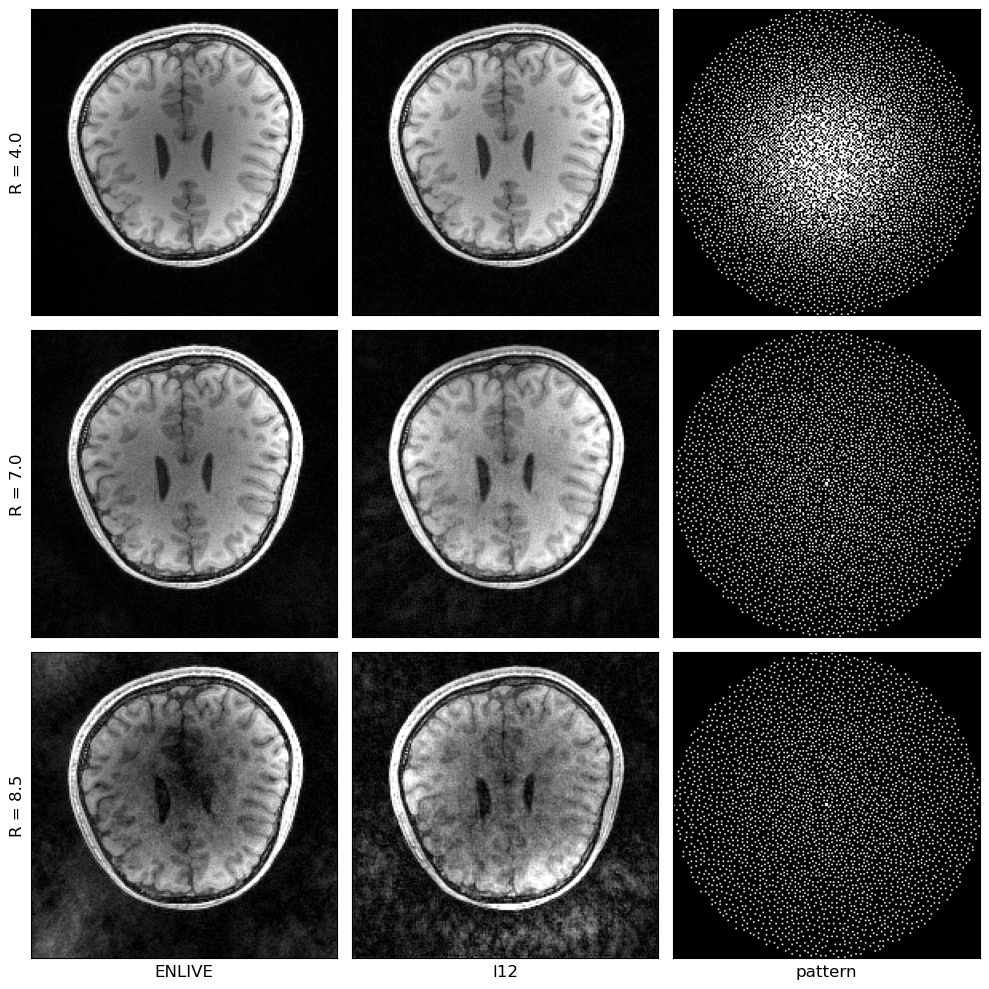}
    \caption{Variable-density Poisson-disc undersampled data with varying undersampling factors were reconstructed using ENLIVE (with $k=2$ in \eqref{ENLIVE}) as well as our proposed method. For undersampling factors up to \(R = 7.0\), both methods produce images that are essentially realistic. However, at an undersampling factor of \(R = 8.5\), the ENLIVE reconstruction becomes very noisy, and creates some features void in the center of the image. With higher undersampling factors, artifacts appear in the images reconstructed with our proposed method, but the brain retains more of its features and shows a more uniform contrast.
}
    \label{fig:10}
\end{figure}

\begin{figure}
    \centering
    \includegraphics[width=0.75\linewidth]{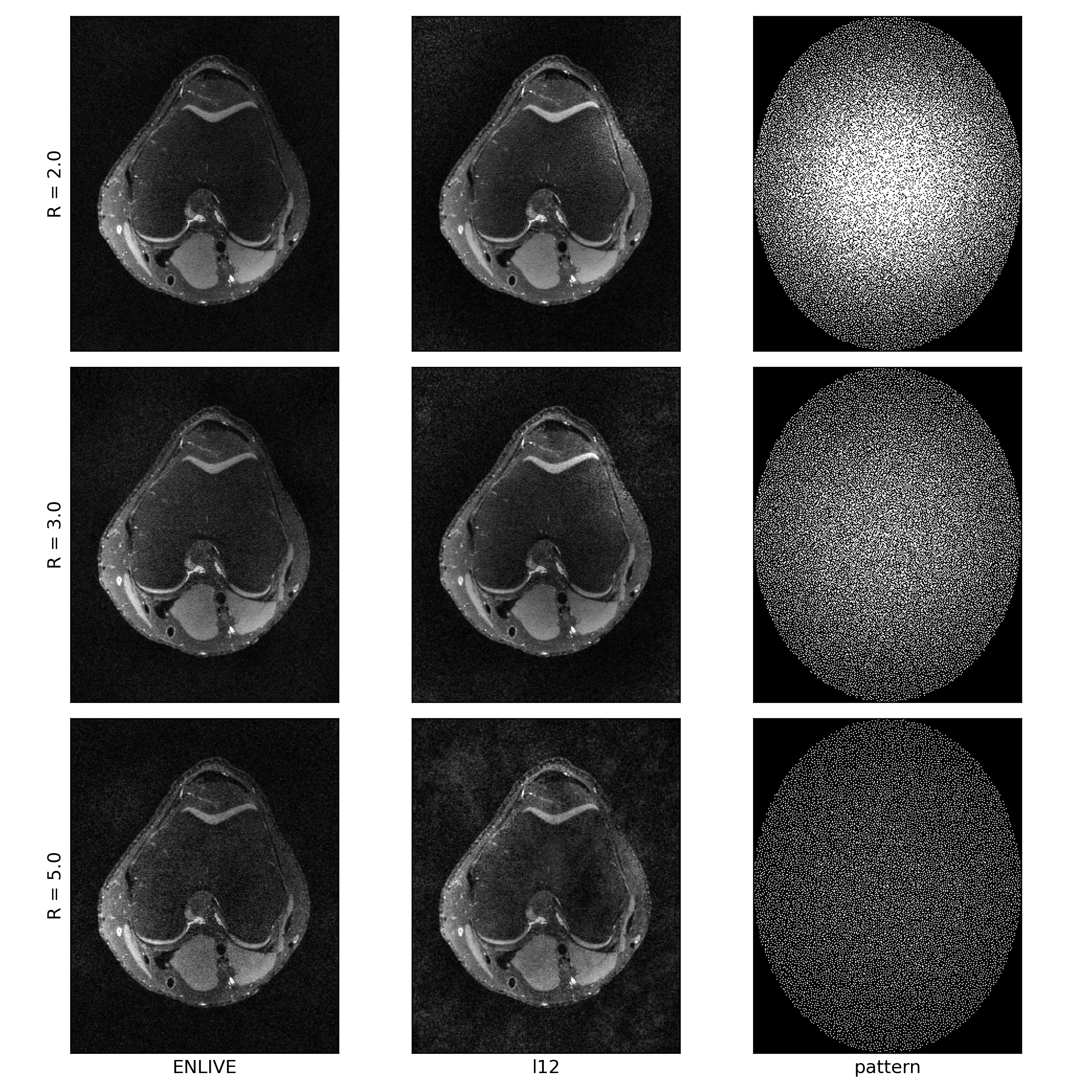}
    \caption{Variable-density Poisson-disc undersampled data of a human knee with varying undersampling factors reconstructed with ENLIVE allowing two sets of maps (i.e. with $k=2$ in eq:\eqref{ENLIVE}) and with our method. The proposed method provides reconstructions with higher and more consistent contrast. For R = 5.0, the reconstruction by ENLIVE seems preferable.
}
    \label{fig:11}
\end{figure}

\begin{figure}
    \centering
    \includegraphics[width=0.75\linewidth]{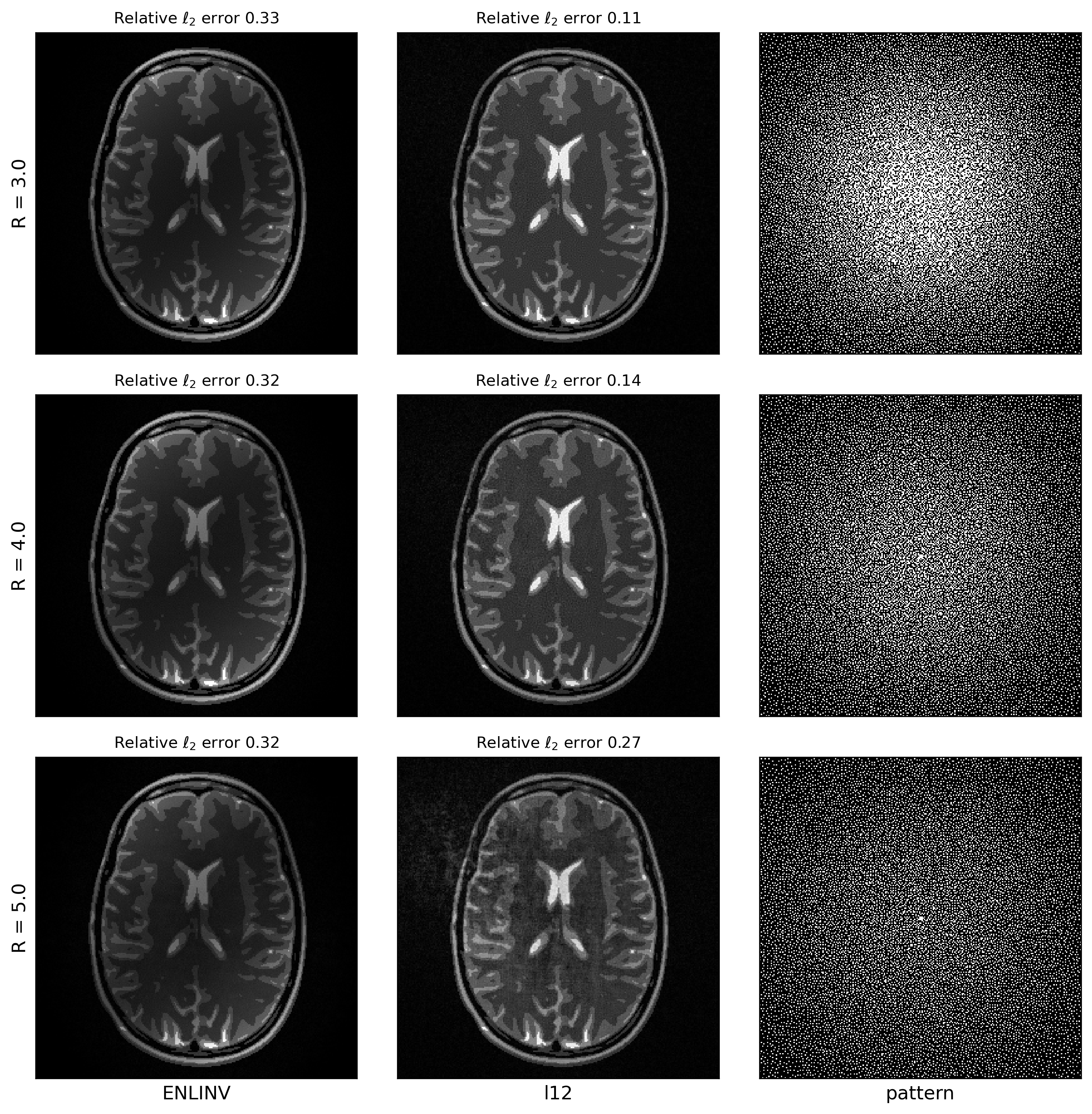}
    \caption{Uniformly undersampled noisy data of the analytical phantom and simulated coils with varying undersampling factors reconstructed with ENLIVE and our method. Relative $\ell_2$ error of the reconstructed phantom and the true phantom is computed. For reduction factors up to $R = 5.0$, the ENLIVE produces reconstructions with low signal at the center, whereas our method produces images with uniform contrast throughout. The reconstructions of our method all have smaller $\ell_2$ error than the ENLIVE reconstructions. For $R=5.0$, artifacts start to appear in our reconstruction.}
    \label{fig:12}
\end{figure}

Next, we compare our reconstruction results with the ENLIVE formulation \eqref{ENLIVE} in the paper by Holme et al. \cite{Holme_Rosenzweig_Ong_Wilke_Lustig_Uecker_2019}. We reproduced the reconstructed brain and knee images exactly as shown in the original paper using the data and source code provided in \url{https://github.com/mrirecon/enlive}. We then compared these results with our own reconstructions. Figures \ref{fig:10} and \ref{fig:11} show calibrationless variable-density Poisson-disc undersampled reconstructions with differing undersampling factors comparing ENLIVE to our proposed method. Figure \ref{fig:12} shows the reconstructions of the simulated analytical phantom from uniformly subsampled Fourier space with different undersampling factors. As a reference, the corresponding sampling patterns are shown in the third columns of each of the plots. 
   
    \newpage
    \section{Proofs}\label{6}
        \subsection{Proof of Lemma \ref{lemma4}}
Recall that $\boldsymbol{\Phi}  = \begin{bmatrix} (\bf{\bar{q}}_1 \bigotimes \bf{b}_1)^* \\ \vdots \\  (\bf{\bar{q}}_N \bigotimes \bf{b}_N)^* \end{bmatrix}$, where each row of $\boldsymbol{\Phi}$ depends only on one row of $\bb{B}$ and $\boldsymbol{\Psi}^*$. Due to symmetry, the probability distribution of $\bold{\Phi}$ is unchanged if we assume $\bb{B}$ is fixed but the rows of $\mathbf{\Psi}^*$ are random. More precisely, we assume the rows  of $\mathbf{\Psi}^* \in \mathbb{C}^{N\times N}$ are chosen uniformly at random without replacement from some $\mathbf{\Psi}_0^* \in \mathbb{C}^{N\times N}$ with $\mathbf{\Psi}_0^* \mathbf{\Psi}_0 = \bb{I}_N$. Under such a model, the rows of $\mathbf{\Psi}^*$ are not independent, making it difficult to analyze the probability model on $\boldsymbol{\Phi}$ directly. On the other hand, if the rows of $\boldsymbol{\Psi}$ are chosen uniformly at random with replacement from $\boldsymbol{\Psi}_0$, we can apply a matrix Bernstein inequality to estimate the operator norm of $\boldsymbol{\Psi}_S^*\boldsymbol{\Psi}_S - \bb{I}_S$. We first state the following Lemma which follows almost exactly from Lemma 4.3 in \cite{ling_strohmer_2015}. The only minor change in the proof is in the estimation of the bounds of $R$ and $\delta^2$ for the matrix Bernstein inequality. 
\begin{lemma} \label{repeated}\footnote{In our application, we use $M=N$.}
For any fixed matrix $\bold{B}\in \mathbb{C}^{N\times k}$ and $\boldsymbol{\Psi}^*\in \mathbb{C}^{N \times M} $ with its rows chosen uniformly at random with replacement from the rows of $\boldsymbol{\Psi}^*_0 \in \mathbb{C}^{N \times M} , \boldsymbol{\Psi}_0^* \boldsymbol{\Psi}_0 = \mathbf{I}_N $. For the matrix $\boldsymbol{\Phi}$ defined in \eqref{Phi}, 
   \begin{equation}
       \|{{\mathbf{\mathbf{\Phi}}}_S^*\mathbf{\mathbf{\Phi}}_S - \mathbf{I}_S}\|_{2\rightarrow 2} \leq \delta
   \end{equation} 
   with probability at least $1-N^{-\alpha}$ if
   $N \geq C_\alpha {\mu^B_{max}}^2{\mu^\mathbf{\mathbf{\mathbf{\mathbf{\Psi}}}}_{max}}^2 kn \log N /\delta^2$.
\end{lemma}

\noindent
The following lemma relates Lemma \ref{repeated} to the desired result.

\begin{lemma}{}
    For any integer $m \leq N$, let the random set $T' := \{t_1',\cdots, t_m'\}$, where $t_l' \in [N]$ are selected independently and uniformly at random from $[N]$. Let $T_m$ be the random subset of $[N]$ chosen uniformly at random among all subsets of cardinality $m$.
    \\
    Given a $\boldsymbol{\Psi}^*$ with $\boldsymbol{\Psi}^*\boldsymbol{\Psi} = \bb{I}_N$, $T'$ and $T_m$, define $\mathbf{\Psi}_{T'}^*$ and $\mathbf{\Psi}_{T_m}^*$ as matrices consisting of rows of $\mathbf{\Psi}^*$ indexed by $T'$ and $T_m$ respectively.
    \\
    Given any fixed $\bb{B}$, consider the forward matrix $\mathbf{\mathbf{\Phi}}_{T'}$ and $\mathbf{\mathbf{\Phi}}_{T_m}$ formulated using $\mathbf{\mathbf{\mathbf{\mathbf{\mathbf{\Psi}}}}}_{T'}$ and $\mathbf{\mathbf{\mathbf{\mathbf{\mathbf{\Psi}}}}}_{T_m}$ respectively as defined in \eqref{lifting}.
    \\
    We have that \[\mathbb{P}(\lambda_{min} (\mathbf{\mathbf{\Phi}}_{T'}^*\mathbf{\mathbf{\Phi}}_{T'}) =0) \geq \mathbb{P}(\lambda_{min} (\mathbf{\mathbf{\Phi}}_{T_m}^*\mathbf{\mathbf{\Phi}}_{T_m}) =0 ), \]
    where $\lambda_{\min} (\mathbf{\mathbf{\Phi}}_{T_m}^*\mathbf{\mathbf{\Phi}}_{T_m})$ is the smallest eigenvalue of the positive semi-definite matrix.
    \begin{proof}
            \begin{align*}
            &    \mathbb{P}(\lambda_{min} (\mathbf{\mathbf{\Phi}}_{T'}^*\mathbf{\mathbf{\Phi}}_{T'}) =0) 
            \\
            & = \Sigma_{k=1}^{m} \mathbb{P}(\lambda_{min}(\mathbf{\mathbf{\Phi}}_{T'}^*\mathbf{\mathbf{\Phi}}_{T'})=0 | card(T') =k ) \mathbb{P}(card(T')=k)
            \\
            & =  \Sigma_{k=1}^{m} \mathbb{P}(\lambda_{min}(\mathbf{\mathbf{\Phi}}_{T_k}^*\mathbf{\mathbf{\Phi}}_{T_k})=0 ) \mathbb{P}(card(T')=k)
            \\
            &  \geq \mathbb{P}(\lambda_{min}(\mathbf{\mathbf{\Phi}}_{T_m}^*\mathbf{\mathbf{\Phi}}_{T_m})=0 ) \Sigma_{k=1}^{m} \mathbb{P}(card(T')=k)
            \\
            & = \mathbb{P}(\lambda_{min}(\mathbf{\mathbf{\Phi}}_{T_m}^*\mathbf{\mathbf{\Phi}}_{T_m})=0 ) .
            \end{align*}
        The inequality holds because whenever $T \subseteq \hat{T}$, we have $\mathbf{\mathbf{\Phi}}_{T}^*\mathbf{\mathbf{\Phi}}_{T} \preccurlyeq \mathbf{\mathbf{\Phi}}_{\hat{T}}^*\mathbf{\mathbf{\Phi}}_{\hat{T}}$.

    \end{proof}
    
\end{lemma}

\subsection{Proof of Lemma \ref{lemma5}}

Suppose that the rows of $\mathbf{F}_\Omega$ are chosen uniformly at random without replacement from $\mathbf{F}$, where $\mathbf{F^*F}=\mathbf{F}\mathbf{F}^*=\frac{N}{L}$. For a realization of $\mathbf{\mathbf{\Phi}}_S$ with $\|\mathbf{\mathbf{\Phi}}_S^*\mathbf{\mathbf{\Phi}}_S - \mathbf{I}_S\|_{2\rightarrow 2}\leq \eta_1$,
\[\|\mathbf{\mathbf{\Phi}}_S^*\bb{F}_\Omega^*\bb{F}_\Omega\mathbf{\mathbf{\Phi}}_S-\mathbf{\mathbf{\Phi}}_S^*\mathbf{\mathbf{\Phi}}_S\|_{2\rightarrow 2} \leq \eta , \, \text{with probability at least} \, 1-\epsilon \] 
provided, 
\[
L \geq \max \{ C_{4}^2 \frac{(1+\eta_1)}{\eta^2}, C_3 \frac{1}{\eta} \} {\mu_{max}^B}^2 {\mu_{max}^\mathbf{\mathbf{\mathbf{\mathbf{\Psi}}}}}^2 kn\ln(\frac{2C_2^2k^2n^2}{\epsilon})
. \]
\begin{proof}
    Denote the rows of $\bb{F}_\Omega$ as $\bb{a}_l^*$ where $\bb{a}_l^*$ are chosen uniformly at random without replacement from the rows of the N point DFT matrix $\bb{F} $ with $\bb{F}^*\bb{F} = \frac{N}{L}I_N$. Define $\bb{Z}_l := \mathbf{\mathbf{\Phi}}_S^*\bb{a}_l\bb{a}_l^*\mathbf{\mathbf{\Phi}}_S, l =1 \cdots L$, we have $\mathbb{E}(\bb{Z}_l)= \mathbf{\mathbf{\Phi}}_S^*\mathbf{\mathbf{\Phi}}_S / L$ and $\sum\limits_{l=1}^L \bb{Z}_l - \mathbb{E}(\bb{Z}_l) = \mathbf{\mathbf{\Phi}}_S^*\bb{F}_\Omega^*\bb{F}_\Omega\mathbf{\mathbf{\Phi}}_S-\mathbf{\mathbf{\Phi}}_S^*\mathbf{\mathbf{\Phi}}_S$
\\
By symmetrization, for any $p\geq 2$
$$
\mathbb{E} \|\sum\limits_{l=1}^L \bb{Z}_l - \mathbb{E}(\bb{Z}_l)\|_{2\rightarrow 2}^p \leq 2^p \mathbb{E}\|\sum\limits_{l=1}\epsilon_l\bb{Z}_l\|_{2\rightarrow2}^p,
$$
where $\epsilon$ is a Rademacher sequence independent of $Z$.
\\
Next, we use the following tail bound for matrix valued Rademacher sums, 
$$
\mathbb{P}_{\epsilon}(\|\sum\limits_{l=1}^L \epsilon_l\mathbf{\mathbf{\Phi}}_S^* \bb{a}_l\bb{a}_l^*\mathbf{\mathbf{\Phi}}_S\|_{2\rightarrow2} \geq t\sigma) \leq 2kn e^{-t^2/2}, t\geq 0, 
$$
where $\sigma = \|\sum\limits_{l=1}^L \mathbf{\mathbf{\Phi}}_S^* \bb{a}_l\bb{a}_l^*\mathbf{\mathbf{\Phi}}_S \mathbf{\mathbf{\Phi}}_S^* \bb{a}_l\bb{a}_l^*\mathbf{\mathbf{\Phi}}_S\|_{2\rightarrow2}^{1/2} \leq \max\limits_{l}\{|\mathbf{\mathbf{\Phi}}_S^*\bb{a}_l|\}\|\bb{F}_\Omega \mathbf{\mathbf{\Phi}}_S\|_{2\rightarrow 2} = \|\mathbf{\mathbf{\Phi}}_S^*\bb{F}_\Omega^*\|_{1\rightarrow 2} \|\bb{F}_\Omega\mathbf{\mathbf{\Phi}}_S\|_{2\rightarrow 2}$. The derivation of the bound can be found, e.g. in Proposition 8.20 of \cite{foucart_rauhut_2015}.
\\
By Proposition 7.13 in~\cite{foucart_rauhut_2015} regarding the moment estimation,
\[
(\mathbb{E}_{\epsilon} \|\sum\limits_{l=1}^L \epsilon_l \bb{Z}_l\|_{2\rightarrow 2}^p)^{1/p} \leq C_1(C_2kn)^{1/p}\sqrt{p}\|\mathbf{\mathbf{\Phi}}_S^*\bb{F}_\Omega^*\|_{1\rightarrow2}\|\bb{F}_\Omega\mathbf{\mathbf{\Phi}}_S\|_{2\rightarrow2}
,\]
where $C_1 = e^{1/(2e)-1/2} \approx 0.73$, $C_2 = 2\sqrt{\pi} e^{1/6} \approx 4.19$. Hence, by H\"{o}lder's inequality,
\[
\mathbb{E} \|\sum\limits_{l=1}^L \epsilon_l \bb{Z}_l\|_{2\rightarrow 2}^p) \leq C_1^p(C_2kn)p^{p/2}\mathbb{E}(\|\mathbf{\mathbf{\Phi}}_S^*\bb{F}_\Omega^*\|_{1\rightarrow2}^{2p})^{1/2}\mathbb{E}(\|\bb{F}_\Omega\mathbf{\mathbf{\Phi}}_S\|_{2\rightarrow2}^{2p})^{1/2}
\]
Combining the bounds, we have
\[
(\mathbb{E} \|\sum\limits_{l=1}^L \bb{Z}_l - \mathbb{E}(\bb{Z}_l)\|_{2\rightarrow 2}^p)^{1/p} \leq 2 (\mathbb{E}\|\sum\limits_{l=1}\epsilon_l\bb{Z}_l\|_{2\rightarrow2}^p)^{1/p} \leq 2 C_1(C_2kn)^{1/p}\sqrt{p}\mathbb{E}(\|\mathbf{\mathbf{\Phi}}_S^*\bb{F}_\Omega^*\|_{1\rightarrow2}^{2p})^{1/2p}\mathbb{E}(\|\bb{F}_\Omega\mathbf{\mathbf{\Phi}}_S\|_{2\rightarrow2}^{2p})^{1/2p}
.\]
Next, we estimate $\mathbb{E}(\|\mathbf{\mathbf{\Phi}}_S^*\bb{F}_\Omega^*\|_{1\rightarrow2}^{2p})^{1/2p}$ and $\mathbb{E}(\|\bb{F}_\Omega\mathbf{\mathbf{\Phi}}_S\|_{2\rightarrow2}^{2p})^{1/2p}$ respectively. For the first bound, observe that $\|\mathbf{\mathbf{\Phi}}_S^*\bb{F}_\Omega^*\|_{1\rightarrow2} \leq \|\mathbf{\mathbf{\Phi}}_S^*\|_{1\rightarrow 2} \|\bb{F}_\Omega^*\|_{1\rightarrow 1}$. Also, recall that  $\mathbf{\mathbf{\Phi}}_S^* = \begin{bmatrix} \bf{\bar{q}_{1,S}} \bigotimes \bf{b_1} &| \cdots &|  \bf{\bar{q}_{N,S}} \bigotimes \bf{b_N} \end{bmatrix}$. Hence,  $\|\mathbf{\mathbf{\Phi}}_S^*\|_{1\rightarrow 2} = \max\limits_{i}  \|\mathbf{q}_{i,s} \otimes \mathbf{b}_i\|_2 = \max\limits_{i}  \|\bb{q}_{i,s}\|_2 \max\limits_{i}\|\bb{b}_i\|_2$. By definition, $\mu_{max}^{\mathbf{B}} = \sqrt{N}\max\limits_{i,j}|\bb{B}_{i,j}|$ and  $\mu_{max}^\mathbf{\mathbf{\mathbf{\mathbf{\Psi}}}} = \sqrt{N}\max\limits_{i,j}|\mathbf{\mathbf{\mathbf{\mathbf{\Psi}}}}_{i,j}|$, we have $\|\mathbf{\mathbf{\Phi}}_S^*\|_{1\rightarrow 2} \leq \frac{\sqrt{kn}}{N} \mu_{max}^B \mu_{max}^\mathbf{\mathbf{\mathbf{\mathbf{\Psi}}}}$.
\\

\noindent
To bound  $\|\bb{F}_\Omega^*\|_{1\rightarrow 1}$, we explicitly write out that 
$$\|\bb{F}_\Omega^*\|_{1\rightarrow 1} = \max\limits_{\|x\|_1 = 1} \sum\limits_{j \in [N]}|\sum\limits_{i \in [L]} \overline{\bb{F}}_{i,j} x_i|\leq \max\limits_{\|x\|_1 = 1}  N \max\limits_j|\sum \limits_{i}\overline{\bb{F}}_{i,j}x_i|\leq N \max\limits_{i,j}|\bb{F}_{i,j}| \leq N/\sqrt{L},$$ since the rows of $\bb{F}_\Omega$ are chosen uniformly from the N point DFT matrix $\bb{F}$ for which $\bb{F}^*\bb{F} = \frac{N}{L} I$.
\\

\noindent
For the term $\mathbb{E}(\|\bb{F}_\Omega\mathbf{\mathbf{\Phi}}_S\|_{2\rightarrow2}^{2p})^{1/2p}$, we first derive a Bernoulli version where each rows of $\bb{F}_\Omega$ are the rows of $\bb{F}$ multiplied by independent Bernoulli selectors. More precisely, consider  $\bb{F}_\Omega=\bb{P}\bb{F}$, where $\bb{P}=\text{diag}(\delta_1,\cdots,\delta_N)$ and $\delta_i$ are independent Bernoulli variables with  $\mathbb{E}(\delta_i) = L/N$. 
\\

\noindent
Apply Lemma 14.3 \cite{foucart_rauhut_2015}  on the non-square matrix $\mathbf{\mathbf{\Phi}}_S^*\bb{F}_\Omega^*$, for any $p\geq 1$,
\[
\mathbb{E}_\delta(\|\mathbf{\mathbf{\Phi}}_S^*\bb{F}^*\bb{P}\|_{2\rightarrow2}^{2p})^{1/2p}\leq \sqrt{2}C_1(C_2kn)^{1/p}\sqrt{p}\mathbb{E}(\|\mathbf{\mathbf{\Phi}}_S^*\bb{F}_\Omega^*\|_{1\rightarrow2}^{2p})^{1/2p} +\sqrt{\frac{L}{N}}\|\mathbf{\Phi}_S^*\bb{F}_\Omega^*\|_{2\rightarrow2}
.\]
Also since $\|\bb{F}_\Omega\mathbf{\Phi}_S\|_{2\rightarrow2} \leq \|\bb{F}_\Omega\|_{2\rightarrow 2}\|\mathbf{\Phi}_S\|_{2\rightarrow2} \leq \sqrt{\frac{N}{L}}\|\mathbf{\Phi}_S\|_{2\rightarrow2}$.Combining the bounds, we have
\[
\begin{split}
&(\mathbb{E} \|\sum\limits_{l=1}^L \bb{Z}_l - \mathbb{E}(\bb{Z}_l)\|_{2\rightarrow 2}^p)^{1/p} 
\\
& \leq
2C_1 (C_2kn)^{1/p}\sqrt{p} \sqrt{\frac{kn}{L}}\mu_{max}^B\mu_{max}^\mathbf{\mathbf{\mathbf{\mathbf{\mathbf{\Psi}}}}} (\sqrt{2}C_1 (C_2kn)^{1/p}\sqrt{p}\sqrt{\frac{kn}{L}}\mu_{max}^B\mu_{max}^\mathbf{\mathbf{\mathbf{\mathbf{\mathbf{\Psi}}}}} + \|\mathbf{\Phi}_S\|_{2\rightarrow 2})
\\
&\leq
(C_2kn)^{2/p}(2\sqrt{2}C_1^2 \frac{knp}{L} {\mu_{max}^B}^2{\mu_{max}^\mathbf{\mathbf{\mathbf{\mathbf{\Psi}}}}}^2 +2C_1\sqrt{\frac{knp}{L}}\mu_{max}^B \mu_{max}^\mathbf{\mathbf{\mathbf{\mathbf{\Psi}}}} \|\mathbf{\Phi}_S\|_{2\rightarrow2} )
.\end{split}
\]
\noindent
Then, we apply a tail bound from moment estimation, see Proposition 7.15 from \cite{foucart_rauhut_2015}. For $u\geq 1$,
\[
\mathbb{P}(\|\sum\limits_{l=1}^L \bb{Z}_l-\mathbb{E}(\bb{Z}_l)\| > e(\alpha_1u+\alpha_2 \sqrt{u}))\leq \beta e^ {-u}
, \]
where $\beta = C_2^2k^2n^2$,  $\alpha_1 = 2\sqrt{2}C_1^2 \frac{kn}{L} {\mu_{max}^B}^2{\mu_{max}^\mathbf{\mathbf{\mathbf{\mathbf{\Psi}}}}}^2$, $\alpha_2 = 2C_1\sqrt{\frac{kn}{L}}\mu_{max}^B \mu_{max}^\mathbf{\mathbf{\mathbf{\mathbf{\Psi}}}} \|\mathbf{\Phi}_s\|_{2\rightarrow2}$.
\\

\noindent
The result implies, for any realization of $\mathbf{\Phi}_S$, $\|\mathbf{\Phi}_s^*\bb{F}_\Omega^*\bb{F}_\Omega\mathbf{\Phi}_S-\mathbf{\Phi}_s^*\mathbf{\Phi}_s\|\leq \eta$ with probability at least $1-\epsilon$ provided,
\begin{equation*}
    e\alpha_1\ln(C_2^2k^2 n^2/\epsilon)\leq \eta/2 \text{ and } e\alpha_2\sqrt{\ln(C_2^2k^2 n^2/\epsilon)}\leq \eta/2 
. \end{equation*}
By the definition of $\alpha_1$ and $\alpha_2$,
\[
\begin{split}
 e\alpha_1\ln(C_2^2k^2 n^2/\epsilon)\leq \eta/2 
 \rightarrow & L  \geq \frac{C_3{\mu_{max}^B}^2 {\mu_{max}^\mathbf{\mathbf{\mathbf{\mathbf{\Psi}}}}}^2 kn \ln(\frac{C_2^2k^2n^2}{\epsilon})}{\eta}  , C_3 =4\sqrt{2}eC_1^2.
 \\
 e\alpha_2\sqrt{\ln(C_2^2k^2 n^2/\epsilon)}\leq \eta/2 \rightarrow &\sqrt{\frac{kn}{L}}\|\mathbf{\Phi}_S\|\leq \frac{\eta}{C_{4}\mu_{max}^B \mu_{max}^\mathbf{\mathbf{\mathbf{\mathbf{\Psi}}}}\sqrt{\ln(C_2^2k^2n^2/\epsilon)}}, C_{4} = 4eC_1 
\\
&\rightarrow L \geq \frac{C_{4}^2 {\mu_{max}^B}^2 {\mu_{max}^\mathbf{\mathbf{\mathbf{\mathbf{\Psi}}}}}^2 kn \ln(C_2^2k^2n^2/\epsilon)\|\mathbf{\Phi}_S\|_{2\rightarrow 2}^2 }{\eta^2} .
\end{split}\]
\noindent
Then conditioned on $\|\mathbf{\Phi}_S^*\mathbf{\Phi}_S - \bb{I}_S\|\leq \eta_1$, we have $\|\mathbf{\Phi}_S\|_{2\rightarrow 2}^2 \leq 1+\eta_1$, the two conditions are equivalent to requiring,
\[
L \geq \max \{ C_{4}^2 \frac{(1+\eta_1)}{\eta^2}, C_3 \frac{1}{\eta} \} {\mu_{max}^B}^2 {\mu_{max}^\mathbf{\mathbf{\mathbf{\mathbf{\Psi}}}}}^2 kn\ln(\frac{C_2^2k^2n^2}{\epsilon}).
\]
This completes the proof for the Bernoulli case. For the case when the rows of $\bb{F}_\Omega$ are chosen uniformly from $F$, use a similar argument as in the later part of the proof of Theorem 14.1 \cite{foucart_rauhut_2015}, where it shows the probability for the uniform model is bounded by twice the one for the Bernoulli model derived here. 

\end{proof}

\subsection{Proof of Lemma \ref{lemma6}}
\begin{proof}
    Partition $\bb{A}$ into blocks of one row and $N$ columns. $\bb{A} = [\bb{A}_1, \bb{A}_2,\cdots, \bb{A}_N]$, where by definition $\bb{A}_i = \bb{F}_\Omega\mathbf{\Phi}_i$.
    \[
    \begin{split}
       & \max\limits_{j \in [|S^c|]}\|\bb{A}_S^*\bb{A}_{S^c_j}\|_{2\rightarrow 2}
       \\
       & \leq 
       \|\bb{A}_S\|_{2\rightarrow 2} \max\limits_{i}\|\bb{A}_i\|_{2\rightarrow 2}
       \\
       & \leq \sqrt{1+\eta}\max\limits_{i}\|\bb{A}_i\|_{2\rightarrow 2} .
    \end{split}
    \]
    The last inequality is implied by $\|\mathbf{A}_S^*\mathbf{A}_S - \mathbf{I}_S\|_{2\rightarrow2} \leq\ \eta$.
    \\

    \noindent
    Since $\bb{A}_i = \bb{F}_\Omega \mathbf{\Phi}_i = \bb{F}_\Omega \text{diag}(\mathbf{\mathbf{\mathbf{\mathbf{\Psi}}}}_i) \bb{B}$, we have $\|\bb{A}_i\|_{2\rightarrow2} \leq  \|\bb{F}_\Omega\|_{2\rightarrow2}\|\mathbf{\Phi}_i\|_{2\rightarrow2}$. Under the assumptions, the rows of $\bb{F}_\Omega$ are chosen from the DFT matrix $\bb{F}$ with $\bb{F}^*\bb{F} = \frac{N}{L}\bb{I}_N$, we have $\bb{F}_\Omega \bb{F}_\Omega^* = \frac{N}{L}\bb{I}_L \rightarrow \|\bb{F}_\Omega \|_{2\rightarrow 2} = \sqrt{\frac{N}{L}}$. Recall that $\mu_{max}^\mathbf{\mathbf{\mathbf{\mathbf{\Psi}}}} = \sqrt{N}\max_{i,j}|\mathbf{\mathbf{\mathbf{\mathbf{\Psi}}}}_{i,j}|$, we have for any $i,j$, $|\mathbf{\mathbf{\mathbf{\mathbf{\Psi}}}}_{i,j}|^2 \leq\frac{1}{N} (\mu_{max}^\mathbf{\mathbf{\mathbf{\mathbf{\Psi}}}})^2$. Hence, for any $i$, $\|\mathbf{\Phi}_i\|_{2\rightarrow 2}^2 = \|\mathbf{\Phi}_i^* \mathbf{\Phi}_i\| = \|\bb{B}^* \text{diag}([|\mathbf{\mathbf{\mathbf{\mathbf{\Psi}}}}_{i,1}|^2,\cdots,|\mathbf{\mathbf{\mathbf{\mathbf{\Psi}}}}_{i,N}|^2]) \bb{B}\|_{2\rightarrow 2} 
    \leq \frac{1}{N}(\mu_{max}^\mathbf{\mathbf{\mathbf{\mathbf{\Psi}}}})^2 \|\bb{B}^*\bb{B}\|_{2\rightarrow 2} =  \frac{1}{N}(\mu_{max}^\mathbf{\mathbf{\mathbf{\mathbf{\Psi}}}})^2$. The last equality is due to the assumption that $\bb{B}^*\bb{B} = I_k$. 
    \\

    \noindent
    Combining the bounds,
    \[
    \max\limits_{j \in [|S^c|]}\|\bb{A}_S^*\bb{A}_{S_j}\|_{2\rightarrow 2} \leq \sqrt{\frac{1+\eta}{L}} \mu_{max}^\mathbf{\mathbf{\mathbf{\mathbf{\Psi}}}}.
    \]
    \end{proof}

\subsection{Proof of Lemma \ref{lemma7}}
    \begin{proof}
    Let $\bb{V} = \bb{A}_S(\bb{A}_S^*\bb{A}_S)^{-1}\sgn(\bb{X}_0)$, define 
    \[\boldsymbol{\mathcal{Y}} := \bb{A}^*\bb{A}_S(\bb{A}_S^*\bb{A}_S)^{-1}\sgn(\bb{X}_0) .\]
    \noindent
    $\boldsymbol{\mathcal{Y}}$ is the exact dual certificate since $\mathcal{P}_S \boldsymbol{\mathcal{Y}} =\sgn(\bb{X}_0)$. Next, we have
    \begin{equation}
            \begin{split}
                \|\bb{V} \|_F
                & = \|\bb{A}_S(\bb{A}_S^*\bb{A}_S)^{-1}\sgn(\bb{X}_0) \|_F
                \leq
                \frac{\sqrt{1+\eta}}{1-\eta} \sqrt{n} .
            \end{split}
        \end{equation}
    The inequality comes from the implication by the conditioning of $\bb{A}_S$ that $\bb{A}_S^*\bb{A}_S$ is invertible with $\|(\bb{A}_S^*\bb{A}_S)^{-1}\|_{2\rightarrow2}\leq \frac{1}{1-\eta}$ and $\|\bb{A}_S\|_{2\rightarrow2} \leq \sqrt{1+\eta}$. 
        
    \end{proof}

\subsection{Proof of Lemma \ref{lemma8}}
We consider the case when the spanning coefficients are chosen uniformly at random from a fixed basis. The structure of the proof is similar to the Gaussian case.
Suppose the columns of $\bb{H}$ are chosen uniformly at random with replacement from the columns of $\bb{W} \in \mathbb{C}^{k\times k }$ where $\bb{W}^*\bb{W} = \bb{I}_k$.  

\begin{proof}
    Let $\bb{A}^\dag = (\bb{A}_S^*\bb{A}_S)^{-1}\bb{A}_S^*$. Suppose there exists an $\alpha$ such that $\max\limits_{j \in S^c} \|\bb{A}^\dag \bb{A}_{S_j}\|_{F} \leq \alpha $. We postpone the estimate of the $\alpha$ to the latter part of the proof. 
\\

\noindent
Define \boldmath $\tilde X$\unboldmath $=  \begin{bmatrix}
        \begin{array}{c | c | c | c}
        \mathbf{h_1}\sgn(z_1)   & \mathbf{h_2}\sgn(z_1)   & \cdots & \mathbf{h_c}\sgn(z_1)  
        \\
        \hline
         \mathbf{h_1}\sgn(z_2)   & \mathbf{h_2}\sgn(z_2)   & \cdots & \mathbf{h_c}\sgn(z_2)  
         \\
         \hline
         \vdots & \vdots & \vdots
         \\
          \mathbf{h_1}\sgn(z_N)   & \mathbf{h_2}\sgn(z_N)   & \cdots & \mathbf{h_c}\sgn(z_N)  
        \end{array}
    \end{bmatrix}$.
Then $\sgn(\bb{X}_0) =$  \boldmath $\tilde X$\unboldmath $/ \|\bb{H}\|_F \in \mathbb{C}^{kN \times C} $
    , where $\mathbf{H} = \begin{bmatrix}
        \mathbf{h_1} | \mathbf{h_2} |  \cdots  |  \mathbf{h_C}
    \end{bmatrix}$, $\mathbf{H} \in \mathbb{C}^{k\times C} .$
\\

\noindent
First, we estimate the norm of $\max\limits_{l\in S^c}\|\sgn(\bb{X}_0)^*\bb{A}^\dag \bb{A}_{S^l}\|_F$:
    \begin{center}
        $\max\limits_{l\in S^c}\|\sgn(\bb{X}_0)^*\bb{A}^\dag \bb{A}_{S^l}\|_F \leq \max\limits_{l \in S^c} \|\bb{A}^\dag \bb{A}_{S^l}\|_F \|\sgn(\bb{X}_0)\|_{2\rightarrow 2} \leq 
    \alpha \|$ \boldmath $\tilde X$\unboldmath $\|_{2\rightarrow 2} / \|\bb{H}\|_F$.
    \end{center}
    \noindent 
   By the Cauchy-Schwarz inequality and the sparsity assumption on $\mathbf{z}$, for any $\bb{y} $, $\|$\boldmath $\tilde X$\unboldmath $ \bb{y}\| \leq \|\bb{H}\|_{2\rightarrow2} \|\sgn(\bb{z})\|_2 \|\bb{y}\|_2 \leq \sqrt{n}\|\bb{H}\|_{2\rightarrow2}$.
    \\

    \noindent 
    We want to bound the two terms, $\|\bb{H}\|_F$ and $\|\bb{H}\|_{2\rightarrow 2}$. For $\|\bb{H}\|_F$, since each column of $\bb{H}$ has norm $1$, we have $\|\bb{H}\|_F = \sqrt{C}$. 
    \\
    
    \noindent
    To bound $\|\bb{H}\|_{2\rightarrow2}$, we apply the non-commutative matrix Bernstein inequality (Theorem 4.4 in \cite{ling_strohmer_2015}) to $\bb{H}\bb{H}^* =  \Sigma_{l = 1}^C \bb{h}_l\bb{h}_l^*$. Consider the independent and centered random matrices $\bb{Z}_l = \bb{h}_l\bb{h}_l^* - \frac{1}{k} \bb{I}_k$ with $\mathbb{E}(\bb{Z}_l) = 0$ and dimension $k \times k$,
    
    \begin{equation}
        R = \max_l \|\bb{Z}_l\| = \|\bb{h}_l\bb{h}_l^* - \frac{1}{k}\mathbf{I}\| \leq \max_l \max \big(\|\bb{h}_l\|^2,\frac{1}{k}\big) = 1,
    \end{equation}
    \begin{equation}
        \sigma^2 =  \max \{\| \sum\limits_{l = 1}^C \mathbb{E}(\bb{Z}_l\bb{Z}_l^*)\|,\| \sum\limits_{l = 1}^C \mathbb{E}(\bb{Z}_l^*\bb{Z}_l)\| \} =  \| \sum\limits_{l = 1}^C \mathbb{E}(\bb{h}_l \underbrace{\bb{h}_l^* \bb{h}_l}_{=1} \bb{h}_l^*)-\frac{2}{k} \mathbb{E}(\bb{h}_l \bb{h}_l^*) + \frac{I}{k^2}\| = \| \sum\limits_{l = 1}^C (\frac{I}{k} - \frac{I}{k^2})\| \leq C/k.
    \end{equation}
    The last equality holds because $\mathbf{h}_l^*\mathbf{h}_l = 1$ and $\mathbb{E}(\mathbf{h}_l\mathbf{h}_l^*) = \frac{1}{k}\mathbf{I}_k$ due to the randomness assumption of $\mathbf{H}$.
    Hence, for any $t>0$, 
    \begin{equation}
        \mathbb{P}(\|\bb{H}\bb{H}^*-\frac{C}{k} \bb{I}\|\geq t)\leq 2k \exp\Big(-\frac{t^2/2}{\sigma^2 + Rt/3}\Big)=2k\exp\Big(-\frac{t^2/2}{\frac{C}{k} + t/3}\Big).
    \end{equation}
    Take $t=\frac{C}{2k}$, we have $\|\bb{H}\| \leq \sqrt{3C/2k}$ with probability at least $1- 2k\exp\big(-\frac{3C}{28k}\big)$.
    \\
    
    \noindent
    Lastly, we estimate the $\alpha$ in the bound of $\max\limits_{j \in S^c} \|\bb{A}^\dag \bb{A}_{S_j}\|_{F} \leq \alpha $. First, conditioned on $\|\mathbf{A}^*_S\mathbf{A}_S - \mathbf{I}_S\| \leq \eta <1,$ 
    \[
\begin{split}
    \max\limits_{j \in S^c} \|\bb{A}^\dag \bb{A}_{S_j}\|_F 
    \leq
     \|(\bb{A}_S^*\bb{A}_S)^{-1}\|_{2\rightarrow2}\|\bb{A}_S\|_{2\rightarrow2}\max\limits_{j \in S^c} \|\bb{A}_{S_j}\|_F 
     \leq 
     \frac{\sqrt{1+\eta}}{1-\eta}\max\limits_{j \in S^c} \|\bb{A}_{S_j}\|_F.
\end{split}
\]
Also, since for any $j$, $\|\bb{A}_{S_j}\|_F = \|\bb{F}_\Omega \diag(\boldsymbol{\Psi}_j)\bb{B}\|_F \leq \|\bb{F}_\Omega\|_{2\rightarrow2}\|\diag(\boldsymbol{\Psi}_j)\bb{B}\|_F \leq \sqrt{\frac{N}{L}} \frac{\mu^\Psi_{max}}{\sqrt{N}} \sqrt{k}$,

Combining the above bounds, we have
\[\max\limits_{j \in S^c} \|\bb{A}^\dag \bb{A}_{S_j}\|_F 
\leq 
\frac{\sqrt{1+\eta}}{1-\eta}\sqrt{\frac{k}{L}} \mu^\Psi_{max}.
 \]
\noindent 
Hence, when $L > 6 {\mu^\Psi_{max}}^2 \frac{1+\eta}{(1-\eta)^2 } n$, we have \begin{equation}
    \max\limits_{l\in S^c}\|\sgn(\bb{X}_0)^*\bb{A}^\dag \bb{A}_{S^l}\|_F \leq \frac{\sqrt{1+\eta}}{1-\eta}\sqrt{\frac{1}{L}} \mu^\mathbf{\mathbf{\mathbf{\mathbf{\Psi}}}}_{max} \sqrt{3n/2 } < 1/2
\end{equation}
with probability at least $1- 2k\exp(-\frac{3C}{28k})$.
\\

\noindent
Equivalently, the result can be stated as $\max\limits_{j \in [|S^c|]} \|\mathcal{P}_{S^c_j}\boldsymbol{\mathcal{Y}}\|_F < 1/2$ with probability at least $1-N^{-\alpha}$ if
\begin{equation*}
    C \geq \Tilde{C}_\alpha \log(N) k,
\end{equation*}
where $\Tilde{C}_\alpha$ grows linearly with respect to $\alpha$.
\\

\noindent
The case when the coefficients $\bb{h}_l,l \in [C]$ are chosen independently and uniformly at random from the complex sphere $S_{\mathbb{C}}^{k-1}$ follows an analogous argument by noticing that $\bb{h}_l =\frac{\bb{s}_l}{\|\bb{s}_l\|_2}$ where $\bb{s}_l$ are independent Gaussian vectors drawn from $\frac{1}{\sqrt{2}}\mathcal{N}(0,1) + i \frac{1}{\sqrt{2}}\mathcal{N}(0,1)$.
\end{proof}

\section*{Acknowledgment}

The authors acknowledge support from NSF DMS-2208356, P41EB032840, and NIH R01HL16351. 

\section*{Code Availability}
The source code for this project will be made openly available at \url{https://github.com/mathgirlstarttocode/AutoCalib_BiConvex_pMRI.git} upon the acceptance of this work for publication. 
    \section{Appendix}
        \subsection{Conditioning of random submatrices}\label{Appendix}
Ideally, we would like to consider the setting in which all the matrices in the forward model are fixed and the randomness comes from the support of the signals $\bb{x} \in \mathbb{C}^N$. We consider the setup where the support $S$ selected uniformly at random among all subsets of $[N]$ of cardinality n.
\\

\noindent
The next lemma is related to the conditioning of random block submatrices for any fixed matrix $\boldsymbol{\Phi}$. We omit the proof in the manuscript.
\begin{lemma}
Suppose $\mathbf{\Phi} = \begin{bmatrix} \mathbf{\Phi}_1 & \cdots & \mathbf{\Phi}_N
\end{bmatrix}$, $\mathbf{\Phi}_i \in \mathbb{C}^{N \times k}$ and $\mathbf{\Phi}_i^*\mathbf{\Phi}_i=\mathbf{I}_k$ for all $i$. Let $S$ be a subset of $[N]$ selected at random according to the Bernoulli model with $\mathbb{E}\text{Card(S)}=n$ or uniformly random from all subset of $[N]$ of cardinality $n$. Let $P=\text{diag}[\delta_1,\dots, \delta_N ]$, where $\delta_i$ are i.i. d Bernoulli variables with $\mathbb{E}\delta_i=\frac{n}{N}$. 
\\

\noindent
If there exists a positive constant $c$, and for $\eta$, $\epsilon \in (0,1)$,   
\[
\begin{split}
& \mu_{block}\leq \frac{c\eta}{\ln(kN/\epsilon)}, \mu_{block}=\max \limits_{i\neq k} \|\mathbf{\Phi}_i^*\mathbf{\Phi}_k\|_2,
\\
& \frac{n}{N}\|\mathbf{\Phi}\|_2^2 \leq \frac{c\eta^2}{\ln(kN/\epsilon)}.
\end{split}
\]
Then, with probability at least ($1-\epsilon$),  $\|\mathbf{\Phi}_S^*\mathbf{\Phi}_S - \mathbf{I}_S\|_{2\rightarrow 2} \leq \eta. $  

\end{lemma}

\noindent
\\
Note that the derivation of the Lemma does not depend on any structures of $\Phi$. However, take into account the structure of our specific $\boldsymbol{\Phi}$, defined as $\boldsymbol{\Phi} = \begin{bmatrix}
    \text{diag}(\boldsymbol{\psi}_1)\bb{B} & \cdots & \text{diag}(\boldsymbol{\psi}_N)\bb{B}
\end{bmatrix}$. The condition on $\mu_{block}$ is strict. Note that the derivation of the Lemma does not depend on any structures of $\Phi$. Recall that our specific $\boldsymbol{\Phi}$ is defined as $\boldsymbol{\Phi} = \begin{bmatrix}
    \text{diag}(\boldsymbol{\psi}_1)\bb{B} & \cdots & \text{diag}(\boldsymbol{\psi}_N)\bb{B}
\end{bmatrix}$. When examining the derived conditions for our specific $\boldsymbol{\Phi}$, we find that it does not readily satisfy the necessary requirements; this makes the Lemma less useful in this context.

        \subsection{Algorithm}\label{Appendix A}
Recall from \eqref{l12},
$\|\mathbf{X}\|_{1,2} = \sum\limits_{j \in [N]} \| \mathbf{X}_{T_j \times [C]}\|_F$. Let $f(\mathbf{X}) = \frac{1}{2}\|\bb{Y}-\bb{A}\bb{X}\|_{F}^2$ and $g(\mathbf{X}) = \lambda\|\bb{X}\|_{1,2}$, an basic update of ISTA \cite{Bubeck_2009} for problem \eqref{unconstraint} at iteration $k$ using stepsize $L$ is,
\begin{align}\label{U1}
\mathbf{X}_k = \argmin_{\mathbf{X}} \{g(\mathbf{X}) + \frac{L}{2} \|\mathbf{X}- [\mathbf{X}_{k-1}-\frac{1}{L}\nabla f(\mathbf{X}_{k-1})]\|_F^2\}.
\end{align}
Together with $\nabla f(\mathbf{X})= \mathbf{A}^*(\mathbf{A}\mathbf{X}-\mathbf{Y})$ and define $\mathbf{Z}_{k-1} = \mathbf{X}_{k-1}-\frac{1}{L}\mathbf{A}^*(\mathbf{A}\mathbf{X}_{k-1}-\mathbf{Y})$, the update \eqref{U1} becomes,
\begin{align}\label{U2}
    \mathbf{X}_k = \argmin_{\mathbf{X}} \{g(\mathbf{X}) + \frac{L}{2} \|\mathbf{X}- \mathbf{Z}_{k-1}\|_F^2\}.
\end{align}
Now, consider the general problem of 
\begin{align}\label{U3}
    \mathbf{X} = \argmin_{\mathbf{X}} \{\lambda\|\mathbf{X}\|_{1,2} + \frac{L}{2} \|\mathbf{X}- \mathbf{Z}\|_F^2\}.
\end{align}
Note that the optimization problem \eqref{U3} is separable. In particular, if we define $\mathbf{X} = \begin{bmatrix}
    \mathbf{X}_1 \\ \vdots \\ \mathbf{X}_N
\end{bmatrix}$, $\mathbf{Z} = \begin{bmatrix}
    \mathbf{Z}_1 \\ \vdots \\ \mathbf{Z}_N
\end{bmatrix}$,
$
    \argmin\limits_{\mathbf{X}} \{\lambda\|\mathbf{X}\|_{1,2} + \frac{L}{2} \|\mathbf{X}- \mathbf{Z}\|_F^2\} = \argmin\limits_{\mathbf{X}_1,\cdots, \mathbf{X}_N} \{\sum\limits_{i \in [N]} \lambda \|\mathbf{X}_i\|_F + \frac{L}{2} \|\mathbf{X}_i - \mathbf{Z}_i\|_F^2 \}
$.
Furthermore, 
\begin{equation}
    \begin{split}\label{update}
    &\argmin_{\mathbf{X}_i} \{\sum\limits_{i \in [N]} \lambda \|\mathbf{X}_i\|_F + \frac{L}{2} \|\mathbf{X}_i - \mathbf{Z}_i\|_F^2 \}
    \\
    = 
    &\argmin_{\mathbf{X}_i} \{\sum\limits_{i \in [N]} \frac{\lambda}{L}\|\mathbf{X}_i\|_F + \frac{1}{2} \|\mathbf{X}_i - \mathbf{Z}_i\|_F^2 \} 
    \\
    =
    &
    \text{prox}_{\frac{\lambda}{L}\|\cdot\|_F}(\mathbf{X}_i) = \mathbf{Z}_i(1-\frac{\lambda /L}{\max(\|\mathbf{Z}_i\|_F,\lambda/L)})
\end{split}
\end{equation}

    \printbibliography
   

\end{document}